\documentclass[pdflatex,sn-mathphys-num]{sn-jnl}% Math and Physical Sciences Numbered Reference Style
\usepackage{graphicx}%
\usepackage{multirow}%
\usepackage{amsmath,amssymb,amsfonts}%
\usepackage{amsthm}%
\usepackage{empheq}
\usepackage{mathtools}
\usepackage{mathrsfs}%
\usepackage[title]{appendix}%
\usepackage{xcolor}%
\usepackage{textcomp}%
\usepackage{manyfoot}%
\usepackage{booktabs}%
\usepackage{algorithm}%
\usepackage{algorithmicx}%
\usepackage{algpseudocode}%
\usepackage{mymacro}

\usepackage{accents}
\usepackage{tikz}
\usetikzlibrary{tikzmark}

\usepackage{aliascnt}
\usepackage[capitalize,noabbrev]{cleveref}
\usepackage{autonum}

\usepackage[sort&compress,numbers]{natbib}

\usepackage{comment} % by fujimoto

\usepackage{xspace}

\makeatletter
\def\bpublisher#1{%
  #1\@ifnextchar,{\SN@publishercomma}{\xspace}%
}
\def\SN@publishercomma,{%
  \@ifnextchar\blocation{\SN@droplocation}{,\space}%
}
\def\SN@droplocation\blocation#1{\xspace}
\makeatother
\theoremstyle{thmstyleone}%
\newtheorem{theorem}{Theorem}%  meant for continuous numbers
\newtheorem{proposition}[theorem]{Proposition}% 
\newtheorem{lemma}[theorem]{Lemma}% 
\newtheorem*{claim}{Claim}% 
\newtheorem{assumption}{Assumption}% 

\theoremstyle{thmstyletwo}%
\newtheorem{remark}{Remark}%

\theoremstyle{thmstylethree}%
\newtheorem{definition}{Definition}%

\crefname{assumption}{Assumption}{Assumptions}
\Crefname{assumption}{Assumption}{Assumptions}

\begin{document}

\title[Wasserstein accelerated proximal gradient algorithms]{Wasserstein accelerated proximal gradient algorithms}


\author*[1]{\fnm{Kansei} \sur{Ushiyama}}\email{ushiyama.k.1053@m.isct.ac.jp}

\author[2]{\fnm{Yi} \sur{Feng}}\email{yi.feng@cs.au.dk}

\author[3]{\fnm{Kaito} \sur{Fujii}}\email{fujii@i.kyoto-u.ac.jp}

\author[4]{\fnm{Yuma} \sur{Fujimoto}}\email{fujimoto\_yuma@cyberagent.co.jp}

\author[5]{\fnm{Noboru} \sur{Isobe}}\email{noboru.isobe@riken.jp}

\affil*[1]{\orgdiv{Department of Mathematical and Computing Science, School of Computing}, \orgname{Institute of Science Tokyo (former Tokyo Institute of Technology)}, \orgaddress{\street{2-12-1 Ookayama}, \city{Meguro-ku}, \postcode{152-8550}, \state{Tokyo}, \country{Japan}}}

\affil[2]{\orgdiv{Department of Computer Science}, \orgname{Aarhus University}, \orgaddress{\street{Åbogade 34}, \city{Aarhus}, \postcode{8200}, \country{Denmark}}}

\affil[3]{\orgdiv{Graduate School of Informatics}, \orgname{Kyoto University}, \orgaddress{\street{Yoshida-honmachi, Sakyo-Ku}, \city{Kyoto-shi}, \postcode{606-8501}, \state{Kyoto}, \country{Japan}}}

\affil[4]{\orgdiv{AI Lab}, \orgname{CyberAgent}, \orgaddress{\street{Udagawa-cho 40-1}, \city{Shibuya-Ku}, \postcode{150-0042}, \state{Tokyo}, \country{Japan}}}

\affil[5]{\orgdiv{Center for Advanced Intelligence Project}, \orgname{RIKEN}, \orgaddress{\street{1-4-1 Nihonbashi}, \city{Chuo-ku}, \postcode{103-0027}, \state{Tokyo}, \country{Japan}}}


\abstract{We propose accelerated proximal gradient methods for composite optimization problems on the quadratic Wasserstein space of probability measures, including Kullback--Leibler divergence minimization.
The objective function consists of a smooth geodesically convex functional and a possibly nonsmooth regularizer that is convex along generalized geodesics.
The convexity and smoothness conditions on the smooth component are imposed only on the feasible set.
As a preliminary result, we establish the equivalence between geodesic convexity of the smooth functional on the feasible set and convexity of its Hilbert-space lift on the corresponding feasible set.
Assuming exact proximal steps, we prove an $\Order{k^{-2}}$ convergence rate for the objective gap when the smooth component is geodesically convex and accelerated linear convergence when it is geodesically strongly convex.
We also describe particle approximations of these algorithms.}




\maketitle

% Manuscript fragment. Requires amsmath and amssymb.
% Bibliography: wasserstein_fb_references.bib
% Drafting scope: the contributions follow the previously discussed
% constrained-proximal construction and O(k^{-2}) objective-gap result.
% Match the regularity assumptions and admissible map class to the final theorem.
% Entropy is used as motivation; applicability requires checking those assumptions.

\section{Introduction}
\label{sec:intro}

We consider composite optimization over the Wasserstein space of probability measures.
Let $\PP_2\paren{\RR^d}$ denote the set of Borel probability measures with finite second moment, equipped with the quadratic Wasserstein distance $W_2$. Let the functionals $\FF$, $\HHr \colon \PP_2(\RR^d) \to \RR\cup\setE{+\infty}$ be the smooth component of the objective and a possibly nonsmooth regularizer, respectively. Let $\iota_{\CC}$ be the indicator of a feasible set $\CC\subseteq\PP_2\paren{\RR^d}$, and set $\HH \coloneqq\HHr+\iota_{\CC}$.
The considered problem is 
\begin{equation}
    \min_{\mu\in\PP_2\paren{\RR^d}} \GG\paren{\mu}\coloneqq\FF(\mu)+\HH(\mu) = \FF(\mu)+\HHr(\mu) + \iota_{\CC}(\mu).
    \label{prob}
\end{equation}
We assume that this problem admits a minimizer $\mu^\star$.

Optimization over probability measures provides a framework for free-energy minimization, sampling, and learning models parameterized by distributions~\cite{JKO98,SKL20,CB18}. A representative example is Kullback--Leibler divergence minimization.
For a target probability measure $\pi=f\dd x \in \PP_2(\RR^d)$ with $f(x)\propto\mathrm{e}^{-V(x)}$, we have
\[
    \argmin_{\mu \in \PP_2(\RR^d)}\mathrm{KL}(\mu\|\pi) = \argmin_{\mu \in \PP_2(\RR^d)}\paren*{\int V \dd \mu + \ent(\mu)},
\]
where $\ent(\mu)\coloneqq\int p(x)\log p(x)\dd x$ for $\mu=p\dd x$, and $\ent(\mu)\coloneqq+\infty$ otherwise.
This problem is an instance of~\eqref{prob} with the potential energy $\FF(\mu)=\int V(x)\dd\mu(x)$, the entropy regularizer $\HHr(\mu)=\ent(\mu)$, and $\CC=\PP_2(\RR^d)$.

The regularization term $\HH$ in~\eqref{prob} need not be smooth or finite everywhere. Entropy, for example, takes the value $+\infty$ on singular measures, and the indicator of a constraint set is infinite outside the feasible region. Thus, a gradient step for the full objective is not generally available. Forward-backward splitting addresses this difficulty by treating the smooth and nonsmooth components separately. In the Wasserstein setting, Salim, Korba, and Luise~\cite{SKL20} proposed a proximal-gradient method that combines a forward step for a smooth convex potential energy with a Wasserstein proximal step for a regularizer convex along generalized geodesics. They established an $\Order{k^{-1}}$ convergence rate for the objective gap, matching the rate of the nonaccelerated proximal-gradient method in Euclidean spaces. In Euclidean optimization, the fast iterative shrinkage-thresholding algorithm (FISTA) improves this rate to $\Order{k^{-2}}$ for convex composite objectives~\cite{BT09}. This motivates the development of an accelerated method for the Wasserstein composite problem.

In this paper, we develop Wasserstein FISTA and establish an $\Order{k^{-2}}$ convergence rate for the objective gap. A feature of our algorithm is that the convexity and smoothness conditions on $\FF$ are imposed only on the feasible set $\CC$, rather than on the entire Wasserstein space. The scope of the proposed algorithm therefore includes functionals that are not geodesically convex on the whole of $\PP_2\paren{\RR^d}$ but are geodesically convex on $\CC$. This distinction is useful when the feasible region is determined by the application and the behavior of the objective outside that region is irrelevant to the optimization problem.

\subsection{Contributions}
\label{subsec:contributions}
We propose two accelerated proximal gradient algorithms on Wasserstein spaces. Their transport-map formulations are given in~\eqref{FISTA} and~\eqref{SR2FISTA} below, and their particle approximations are presented in \cref{alg:FISTA,alg:SR2FISTA}.
\cref{thm:FISTA,thm:SR2FISTA} show the convergence rates for the updates~\eqref{FISTA} and~\eqref{SR2FISTA}. 
\begin{itemize}
    \item For geodesically convex $\FF$, both updates achieve  
    \[
        \GG(\mu_k) - \GG(\mu^\star) \le \frac{2L}{k^2} W_2(\mu_0,\mu^\star)^2.
    \]
    % This algorithm can be seen as an Wasserstein adaptation of an proximal accelerated gradient algorithm in \cite[Algorithm~20]{ADA21}.
    \item For $m$-geodesically convex $\FF$ with $m>0$, they achieve
    \[
        \GG(\mu_k) - \GG(\mu^\star) \le \Order{\paren*{1-\sqrt{\frac mL}}^k} \qquad \GG(\mu_k) - \GG(\mu^\star) \le \Order{\paren*{1+\sqrt{2\frac mL}}^{-k}},
    \]
    respectively.
    % This algorithm can be seen as an Wasserstein adaptation of a variant of SR2FISTA~\cite[Algorithm~1]{U25}.
\end{itemize}
% Wasserstein SR2FISTA has a better convergence rate, whereas Wasserstein FISTA exhibits better numerical performance for some problems.

Moreover, as a preliminary result of the convergence analysis, we show the equivalence between the geodesic convexity of $\FF$ on the feasible set and the convexity of its Hilbert-space lift on the corresponding feasible set in \cref{prop:convex}.

\subsection{Challenges and our approaches}
In this section, we informally introduce challenges in developing accelerated proximal gradient methods on the Wasserstein space and our approaches to addressing them.
A formal treatment is provided in~\cref{sec:notation} and the subsequent sections.
In this subsection, we consider the unconstrained case $\CC = \PP_2(\RR^d)$ for simplicity.

The nonsmooth component $\HH$ presents a central difficulty in extending acceleration in Euclidean space to Wasserstein space.
Using the Hilbert-space lifting in Lions' calculus \cite{C13,CD18}, we fix a reference measure $\nu$ and represent $\mu=X_\#\nu$ with $X\in L^2(\nu)$, introducing the lifted functionals $\tF(X)\coloneqq\FF(X_\#\nu)$ and $\tH(X)\coloneqq\HH(X_\#\nu)$. In the proximal optimization setting where $\FF$ is smooth and geodesic convex, $\tF$ becomes convex on $L^2(\nu)$. If $\tH$ is also convex on $L^2(\nu)$, a direct application of Euclidean FISTA to $\min_{X\in L^2(\nu)}\tF(X)+\tH(X)$ is possible, since its convergence analysis does not require finite dimensionality. However, we assume that $\HHr$ is convex along generalized geodesics, and this does not, in general, imply convexity of $\tH$ on the whole of $L^2(\nu)$. Consequently, a direct extension of Euclidean FISTA does not justify its usual convergence analysis. 
This difficulty is particularly relevant to entropy regularization.
Indeed, existing theoretical results on discrete-time acceleration methods for Wasserstein optimization do not apply to such problems.
Tang, Li, and Gangbo~\cite[Remark~2.4]{TLG26} explicitly note that their framework does not directly cover entropy-type functionals, while Tanaka~\cite[Section~6.2]{T23} does not establish the smoothness condition required for the convergence analysis of entropy-based objectives.

Note also that if both $\tF$ and $\tH$ are convex on $L^2(\nu)$, an optimal measure can be chosen to be a Dirac measure, and the problem reduces to finding its support point using Euclidean optimization algorithms.
This motivates the development of Wasserstein optimization algorithms under the weaker assumption that $\HHr$ is convex along generalized geodesics.

Although $\tH$ is not necessarily convex on $L^2(\nu)$, it is convex on the set $\KK_\nu \subseteq L^2(\nu)$ of optimal transport maps from $\nu$.
We therefore address this difficulty by designing algorithms whose iterates remain in $\KK_\nu$.
We achieve this by restricting the proximal step to $\KK_\nu$ and adapting the constraint-preserving Euclidean algorithm of d'Aspremont, Scieur, and Taylor~\cite[Algorithm~20]{ADA21} to the lifted problem.
Keeping the iterates in $\KK_\nu$ allows us to apply the convexity inequalities for $\tH$ on this set and to analyze convergence using Lyapunov functions similar to those used for the corresponding Euclidean algorithms.
% For maps in $\KK_\nu$, linear interpolation induces a generalized geodesic with base $\nu$. Together with our convexity assumption on the lifted constraint set, 
% The convexity along generalized geodesics of $\HHr$ yields the convexity of the restricted regularizer
% \begin{equation}
%     \Phi\coloneqq\tH+\iota_{\KK_\nu}.
% \end{equation}
% The constructed accelerated formulation updates an auxiliary map through the proximal operator of $\Phi$ and forms the primary iterate and the gradient-evaluation point through convex combinations. Under our assumptions on the lifted constraint set, these updates preserve feasibility and membership in $\KK_\nu$. The resulting method can therefore use the Hilbert-space inequalities for the smooth term together with the convexity inequalities supplied by generalized geodesics for the regularizer. This combination closes the accelerated convergence argument while requiring the convexity and smoothness conditions on $\FF$ only on $\CC$.

\subsection{Related work}
\label{sec:rel}

\subsubsection{Wasserstein optimization}
Wasserstein optimization builds on the JKO (Jordan--Kinderlehrer--Otto) scheme~\cite{JKO98} and the theory of gradient flows in spaces of probability measures~\cite{AGS08}.
For composite objectives, Salim, Korba, and Luise~\cite{SKL20} propose a forward--backward method for a smooth convex potential energy and a regularizer convex along generalized geodesics, establishing an $\Order{k^{-1}}$ objective-gap rate.
Related methods include proximal coordinate algorithms~\cite{YCY24}, forward--backward Gaussian variational inference~\cite{DBCS23}, and Wasserstein mirror and preconditioned gradient methods~\cite{BUDAK24}.
Bonet et al.~\cite[Appendix F]{BUDAK24} also provide a composite Bregman proximal-gradient extension, with regularity conditions tied to the transport interpolations used in the analysis.
For KL-regularized objectives, Lascu et al.~\cite{LMSS26} establish linear convergence of proximal descent schemes under flat convexity and logarithmic Sobolev assumptions.

Acceleration has been studied through continuous-time dynamics~\cite{WL22,CLTW25}, including $\Order{t^{-2}}$ objective-gap guarantees under geodesic convexity and suitable regularity assumptions.
Discrete $\Order{k^{-2}}$ rates have also been established by Tanaka~\cite{T23} using transport-map convexity and smoothness, and by Tang, Li, and Gangbo~\cite{TLG26} through the lifting for smooth, geodesically convex functionals.
% \begin{comment}
    {
    % \color{blue}memo:
    Park et al.~\cite{NEURIPS2025_347dae37} established an $\Order{k^{-\log_2(1+\sqrt{2})}}$ rate for smooth potential functionals in Wasserstein gradient descent.
    }
% \end{comment}

\subsubsection{Accelerated forward--backward methods in Euclidean spaces}

Forward--backward splitting~\cite{P79,RB16} is a classical method for composite convex optimization in Euclidean spaces, where the objective $f=g+h$ consists of a smooth convex function $g\colon\RR^d \to \RR$ with an $L$-Lipschitz continuous gradient and a possibly nonsmooth convex regularizer $h\colon\RR^d \to \RR\cup\setE{+\infty}$. The proximal-gradient method, known as ISTA~\cite{DDD04}, achieves an $\Order{k^{-1}}$ objective-gap rate. FISTA~\cite{BT09} and other accelerated proximal-gradient methods~\cite{N13} accelerate this to $\Order{k^{-2}}$, while OptISTA~\cite{JGR25} attains the optimal worst-case objective-gap bound.

When $g$ is additionally $m$-strongly convex for some $m>0$, accelerated proximal-gradient methods achieve linear convergence with objective-gap rates of the form $\Order{(1-\sqrt{m/L})^k}$~\cite{CP16,CC19,FV19,FV20,ADA21,RC22}. Aujol et al.~\cite{ADR22} improved the leading exponential rate by a factor of $\sqrt{2}$. Subsequently, Ushiyama~\cite{U25} proposed SR2FISTA, which further improves the linear convergence factor while retaining an $\Order{k^{-2}}$ objective-gap guarantee. More recently, Prox-ITEM~\cite{UTDT26} achieves the optimal squared-distance convergence rate $\norm{x_k-x^\star}^2=\Order{(1-\sqrt{m/L})^{2k}}$. However, this distance bound does not directly imply a corresponding objective-gap guarantee in the presence of nonsmooth regularization.

\section{Notation and preliminaries}\label{sec:notation}
\subsection{General notation}
For $x,y \in \RR^d$, $\inpr{x}{y}$ and $\norm{x}$ denote the Euclidean inner product and norm, respectively.

Let $\PP_2(\RR^d)$ denote the space of Borel probability measures on $\RR^d$ with finite second moment, and $\PP_{2}^{\mathrm{ac}}(\RR^d)$ denote its subset of measures that are absolutely continuous with respect to Lebesgue measure.
For a measurable map $T \colon \RR^d \to \RR^n$ and a measure $\mu \in \PP_2(\RR^d)$, the pushforward measure of $\mu$ by $T$ is denoted by $T_\# \mu$, which is defined by
$T_\# \mu (A) \coloneqq \mu(T^{-1}(A))$ for any Borel set $A \subseteq \RR^n$,
where $T^{-1}(A) = \setI{x \in \RR^d}{T(x) \in A}$.
For a functional $\GG \colon \PP_2(\RR^d) \to \RR \cup \setE{+\infty}$, the effective domain of $\GG$ is defined by 
$
    D(\GG) \coloneqq \setI{\mu \in \PP_2(\RR^d)}{\GG(\mu) < +\infty}.
$

For any $\mu \in \PP_2(\RR^d)$, $L^2(\mu)$ denotes the space of functions $T \colon \RR^d \to \RR^d$ such that $\int_{\RR^d} \norm{T(x)}^2 \dd \mu(x) < +\infty$, equipped with the inner product $\inprL{S}{T}{\mu} \coloneqq \int_{\RR^d} \inpr{S(x)}{T(x)}\dd \mu(x)$ and the norm $\normL{T}{\mu} \coloneqq \sqrt{\inprL{T}{T}{\mu}}$.
For a Fr\'echet differentiable functional $\tilde\FF \colon L^2(\mu) \to \RR$, $\nabla \tilde\FF$ denotes its Fr\'echet gradient. For $A \subseteq L^2(\nu)$, $\iota_{A} \colon L^2(\nu) \to \RR\cup\setE{+\infty}$ denotes the indicator function defined by $\iota_A(T) = 0$ if $T \in A$ and $\iota_A(T) = +\infty$ if $T \not\in A$.

For $\mu, \rho \in \PP_2(\RR^d)$, the set of couplings $\Gamma(\mu,\rho)$ denotes the set of $\gamma \in \PP_2(\RR^d\times\RR^d)$ with first marginal $\mu$ and second marginal $\rho$.
The 2-Wasserstein distance $W_2 \colon \PP_2(\RR^d) \times \PP_2(\RR^d) \to \RR$ is defined by
\[
        W_2(\mu,\rho) \coloneqq \inf_{\gamma \in \Gamma(\mu,\rho)} \paren*{\int_{\RR^d \times \RR^d} \norm{x-y}^2 \dd \gamma(x,y)}^{\frac12}.
\]
The set of optimal transport plans $\Gamma_\mathrm{o}(\mu,\rho)$ is defined by
\[
    \Gamma_\mathrm{o}(\mu,\rho) \coloneqq \setI*{\gamma \in \Gamma(\mu,\rho)}{W_2^2(\mu,\rho) = \int_{\RR^d \times \RR^d} \norm{x-y}^2 \dd \gamma(x,y)}.
\]
We define the set of optimal transport maps from $\nu \in \PP_2^\mathrm{ac}(\RR^d)$:
\[
    \KK_\nu \coloneqq \setI{T \in L^2(\nu)}{(\id,T)_\# \nu \in \Gamma_\mathrm{o}(\nu, T_\# \nu)}.
\]

\begin{definition}[Wasserstein gradient]
    A functional $\FF \colon \PP_2(\RR^d) \to \RR\cup\setE{+\infty}$ is said to be \textit{Wasserstein differentiable at} $\mu \in D(\FF)$ if there exists some $\xi \in L^2(\mu)$ such that 
    \begin{equation}
        \FF(\rho) - \FF(\mu) = \int_{\RR^d \times \RR^d} \inpr{\xi(x)}{y-x} \dd \gamma(x,y) + \mathrm{o}(W_2(\mu,\rho)) \quad \text{as} \quad W_2(\mu,\rho)\to0
    \end{equation}
    for all $\rho \in \PP_2(\RR^d)$ and all $\gamma \in \Gamma_\mathrm{o}(\mu,\rho)$.
    For this $\FF$, such $\xi$ of minimal norm is called the \textit{Wasserstein gradient of $\FF$ at $\mu$}, denoted by $\nablaw \FF(\mu)$.
\end{definition}

\bmhead{Lifting to $L^2(\nu)$}
Following~\cite{TLG26}, we consider a measure $\nu \in \PP_{2}^{\mathrm{ac}}(\RR^d)$ as a reference measure and lift $\FF$ and $\HH$ of the problem~\eqref{prob} to functionals $\tF$ and $\tH$ on $L^2(\nu)$, respectively: for $X \in L^2(\nu)$, set $\tHr(X)\coloneqq\HHr(X_\#\nu)$ and define
\begin{equation}
    \tF(X) \coloneqq \FF(X_\# \nu), \quad \tH(X) \coloneqq \tHr(X) + \iota_{\CCnu}(X) = \HHr(X_\# \nu) + \iota_{\CC}(X_\# \nu). \label{tF}
\end{equation}

The differentiability of $\FF$ and $\tF$ is related as follows.
\begin{proposition}\label{prop:diff}
    For $\FF$ and $\tF$ defined by~\eqref{tF}, let $X \in L^2(\nu)$ satisfy $X_\# \nu \in D(\FF)$.
    Then, $\FF$ is Wasserstein differentiable at $X_\# \nu$ if and only if $\tF$ is Fr\'echet differentiable at $X$. Moreover, it holds that
    \[
        \nabla \tF(X) = \nablaw \FF(X_\# \nu)\circ X.
    \]
\end{proposition}
The proof is provided, e.g., in~\cite[Corollary~3.22]{GT19}.

In this paper, we consider the constraint $\CC \subseteq \PP_2(\RR^d)$.
We transfer this constraint into $L^2(\nu)$ by introducing the notation
\[
    \CCnu \coloneqq \setI{X \in L^2(\nu)}{X_\# \nu \in \CC}.
\]

\subsection{Convexity}

For the smooth part $\FF$ of~\eqref{prob}, we impose the smoothness and convexity conditions only on the constraint set $\CC$, as in the Euclidean setting of~\cite[Definition~4.2]{ADA21}.
For convexity on constraints, we introduce the following definition.
\begin{definition}[Geodesic convexity on a subset]\label{def:geoconvex}
    Let $\FF \colon \PP_2(\RR^d) \to \RR\cup\setE{\infty}$. For $m \ge 0$, $\FF$ is said to be \textit{$m$-geodesically convex on} $\CC \subseteq \PP_2(\RR^d)$ if for all $\mu_0, \mu_1 \in \CC$ and $\gamma \in \Gamma_\mathrm{o}(\mu_0,\mu_1)$, it holds that for all $t \in [0,1]$,
    \[
        \FF(\mu_t) \le (1-t)\FF(\mu_0) + t\FF(\mu_1) - \frac m2 t(1-t)W_2^2(\mu_0,\mu_1),
    \]
    where $\mu_t$ is the geodesic defined by $\mu_t \coloneqq {T_t}_\# \gamma$ with $T_t(x_0,x_1) \coloneqq (1-t)x_0 + tx_1$. We simply say geodesically convex when $m=0$.
\end{definition}
\begin{definition}[$L^2$ convexity on a subset]\label{def:L2convex}
    Let $\UU \colon L^2(\rho) \to \RR\cup\setE{+\infty}$ for $\rho \in \PP_2(\RR^d)$. For $m \ge 0$, $\UU$ is said to be \textit{$m$-convex on} $\tilde\CC \subseteq L^2(\rho)$ if for all $T_0, T_1 \in \tilde\CC$, it holds that for all $t \in [0,1]$,
    \[
        \UU(T_t) \le (1-t)\UU(T_0) + t\UU(T_1) - \frac m2 t(1-t)\normL{T_0 - T_1}{\rho}^2,
    \]
    where $T_t$ is defined by $T_t \coloneqq (1-t)T_0 + tT_1$. 
    We simply say convex when $m = 0$.
\end{definition}

For the nonsmooth part $\HH$, we consider convexity along generalized geodesics.
\begin{definition}[Convexity along generalized geodesics]\label{def:genegeoconvex}
    Let $\nu \in \PP_2^{\mathrm{ac}}(\RR^d)$ be a fixed measure. For $\mu_0, \mu_1 \in \PP_2(\RR^d)$, let $T_0,  T_1 \in \KK_\nu$ denote the (unique) optimal transport maps from $\nu$ to $\mu_0, \mu_1$, respectively. \textit{The generalized geodesic from $\mu_0$ to $\mu_1$ with base $\nu$} is defined by
    \[
        \mu_t = ((1-t)T_0 + tT_1)_\#\nu.
    \]
    
    A functional $\HH\colon \PP_2(\RR^d) \to \RR\cup\setE{\infty}$ is said to be \textit{$m$-convex along generalized geodesics with base $\nu$} if it holds that for all $\mu_0, \mu_1$ and $t \in [0,1]$, the generalized geodesic $\mu_t$ from $\mu_0$ to $\mu_1$ with base $\nu$ satisfies
    \[
        \HH(\mu_t) \le (1-t)\HH(\mu_0) + t\HH(\mu_1) - \frac m2 t(1-t)\normL{T_0-T_1}{\nu}^2.
    \]
    We simply say convex along generalized geodesics when $m=0$.
\end{definition}

\section{Assumptions and examples}
\subsection{Assumptions and a preliminary result}
Based on the preceding definition and observation, we consider the following assumptions on the constraint and objective function.
\begin{assumption}[Assumption on constraint]\label{asm:C}
    For the constraint $\CC \subseteq \PP_2(\RR^d)$ of the problem~\eqref{prob}, we assume that $\CCnu$ is closed and convex.
\end{assumption}
\begin{assumption}[Assumption on $\FF$]\label{asm:F}
    For the functional $\FF$ and constraint $\CC$ of the problem~\eqref{prob}, we consider the following conditions.
    \begin{enumerate}
        \item There exists an open neighborhood $\CC' \supseteq \CC$ such that $\CC' \subseteq D(\FF)$ and $\FF$ is Wasserstein differentiable at every $\mu \in \CC'$. \label{asm:diff}
        \item Let $L >0$. For all $S, T \in \CCnu$, it holds that
        \begin{equation}
            \tF(S) - \tF(T) \le \inprL{\nabla \tF(T)}{S - T}{\nu} + \frac L2 \normL{S-T}{\nu}^2. \label{Lsmooth}
        \end{equation}\label{asmLsmooth}
        \item $\FF$ is $m$-geodesically convex on $\CC$. \label{asmconvex}
    \end{enumerate}
\end{assumption}

For \cref{asm:F}.\ref{asmLsmooth}, a sufficient condition for~\eqref{Lsmooth} can be given by the inequality 
\[
    \norm{\nablaw \FF(\mu)(x) - \nablaw \FF(\rho)(y)} \le \frac L2 \paren*{\norm{x-y} + W_2(\mu,\rho)}
\]
for all $\mu,\rho \in \CC$ and $x,y \in \RR^d$. Indeed, for any $S,T\in\CCnu$, we have
\begin{align}
    \MoveEqLeft
    \normL{\nabla \tF(S) - \nabla \tF(T)}{\nu} = \normL{\nablaw \FF(S_\# \nu)\circ S - \nablaw \FF(T_\# \nu)\circ T}{\nu} \\
    &= \paren*{\int \norm{\nablaw \FF(S_\# \nu)(S(x)) - \nablaw \FF(T_\# \nu)(T(x))}^2 \dd \nu(x)}^{\frac12}\\
    &\le \frac{L}{2}\paren*{ \int \paren*{\norm{S(x) - T(x)} + W_2(S_\# \nu, T_\# \nu)}^2 \dd \nu(x)}^{\frac12}\\
    &\le \frac{L}{2}\paren*{\normL{S-T}{\nu} + W_2(S_\# \nu, T_\# \nu)} \le L\normL{S-T}{\nu},
\end{align}
where we used \cref{prop:diff} and $W_2(S_\# \nu, T_\# \nu) \le \normL{S-T}{\nu}$.
The same condition for the case $\CC=\PP_2(\RR^d)$ is provided in~\cite{TLG26}.
Although this condition is intrinsic as it is expressed without using a lift from $\nu$, we instead directly check~\eqref{Lsmooth} for subsequent examples of $\FF$ for convenience.

For \cref{asm:F}.\ref{asmconvex}, the geodesic convexity is a weaker condition than the $L^2(\nu)$-sense convexity in general.
However, they can be equivalent under some conditions.
In~\cite[Theorem~1.1]{P24}, it is shown that if $\FF$ is Wasserstein differentiable, the geodesic convexity of $\FF$ on $\PP_2(\RR^d)$ is equivalent to the convexity of $\tF$ on $L^2(\nu)$.
In this paper, we show that this equivalence still holds even if convexity is restricted to the constraint set under \cref{asm:C}.
\begin{theorem}\label{prop:convex}
    For $\CC \subseteq \PP_2(\RR^d)$, assume \cref{asm:C}.
    For $\FF$ and $\tF$ defined by~\eqref{tF}, assume \cref{asm:F}.\ref{asm:diff}.
    Then, $\FF$ is $m$-geodesically convex on $\CC$ if and only if $\tF$ is $m$-convex on $\CCnu$.
\end{theorem}
The proof is deferred to~\cref{sec:convex}.

Therefore, under~\cref{asm:C}, \ref{asm:F}.\ref{asm:diff}, and \ref{asm:F}.\ref{asmconvex}, it holds that for all $S, T \in \CCnu$,
\begin{equation}
    \tF(S) - \tF(T) \ge \inprL{\nabla \tF(T)}{S - T}{\nu} + \frac m2 \normL{S-T}{\nu}^2. \label{convexF}
\end{equation}

\begin{assumption}[Assumption on $\HH$] \label{asm:H}
    For the functional $\HH = \HHr + \iota_\CC$ of the problem~\eqref{prob}, we assume the following conditions.
    \begin{enumerate}
        \item $D(\HH) \neq \varnothing$.
        \item $\HHr$ is lower semicontinuous.
        \item $\HHr$ is convex along the generalized geodesics with base $\nu$.\label{H3}
    \end{enumerate}
\end{assumption}

From \cref{asm:C} and \ref{asm:H}.\ref{H3}, $\tH$ is convex on $\KK_\nu$.
Thus, we introduce $\Phi \colon L^2(\nu) \to \RR\cup\setE{+\infty}$ such that
\begin{equation}
    \Phi \coloneqq \tH + \iota_{\KK_\nu}, \label{phi}
\end{equation}
so that $\Phi$ is lower semicontinuous and convex on $L^2(\nu)$ since $\KK_\nu$ is closed and convex.
We define its subdifferential by
\[
    \partial \Phi (T) = \setI{G \in L^2(\nu)}{\Phi(S) \ge \Phi(T) + \inprL{G}{S-T}{\nu} \text{ for all } S \in L^2(\nu)}.
\]
Using this subdifferential, it holds that
\begin{equation}
    \tH(S) - \tH(T) \ge \inprL{G}{S - T}{\nu} \quad \text{for any } G \in \partial \Phi(T), \label{convexG}
\end{equation}
for any $S, T \in D(\tH)\cap\KK_\nu$.
This inequality need not hold for $S\notin\KK_\nu$.

% Moreover, we consider a stronger assumption.
% \begin{assumption}[$L^2$ convexity of $\HH$] \label{asm:H2}
%     In addition to \cref{asm:H}, we assume that .
% \end{assumption}
Note that if $\tHr$ is convex on $L^2(\nu)$, including the case $\HH = \iota_{\CC}$, the analogous inequality to~\eqref{convexG} holds for all $S,T\in D(\tH)$ and $G\in\partial\tH(T)$.
Under this stronger assumption, the problem can be reduced to finite-dimensional optimization.
In fact, in this case, the full lifted objective $\tilde\GG=\tF+\tH$ also becomes convex, and for such a problem, a Dirac measure $\delta_{x^\star}$ for some $x^\star \in \RR^d$ can be chosen as an optimal measure $\mu^\star$, and the optimization over $\PP_2(\RR^d)$ can be reduced to finding $x^\star$ of the Dirac measure, which can be solved by Euclidean FISTA.

\subsection{Proximal map}\label{sec:prox}
% Proposed algorithms use two proximal map depending on the convexity of $\HH$.
% For $\HH$ satisfying \cref{asm:H} and \cref{asm:H2}, we use the ordinary proximal map on $L^2(\nu)$:
% \begin{equation}
%     \prox_{\delta\tH}(Y) = \argmin_{X\in L^2(\nu)}\paren*{\tH(X)+\frac{1}{2\delta}\normL{X-Y}{\nu}^2} = \argmin_{X\in\CCnu}\paren*{\tHr(X)+\frac{1}{2\delta}\normL{X-Y}{\nu}^2}.\label{L2prox}
% \end{equation}

In the proposed algorithms, we assume that the proximal map with respect to $\Phi$ is computable: for $\eta > 0$, we define
\begin{align}
    \prox_{\eta\Phi}(Y) &\coloneqq 
    \argmin_{X\in L^2(\nu)}\paren*{\Phi(X) +\frac{1}{2\eta}\normL{X-Y}{\nu}^2}\\
    &= \argmin_{X\in\CCnu\cap\KK_\nu}\paren*{\tHr(X)+\frac{1}{2\eta}\normL{X-Y}{\nu}^2}.\label{geoprox}
\end{align}
Since $\Phi$ is proper, lower semicontinuous, and convex, the proximal map is well-defined.
By Brenier's theorem, $X \in \KK_\nu$ can be expressed as $X = \nabla \phi$ for some convex function $\phi \colon \RR^d \to \RR\cup\setE{+\infty}$.
Therefore,
\begin{align}
    \prox_{\eta\Phi}(Y) = \nabla \phi_Y, \quad \phi_Y \in \argmin_{\substack{\phi \colon \text{convex}\\ \text{s.t. } \nabla \phi \in \CCnu}}\paren*{\tHr(\nabla \phi)+\frac{1}{2\eta}\normL{\nabla \phi-Y}{\nu}^2}. \label{proxB}
\end{align}
This can be approximated by input-convex neural networks (ICNN), as in JKO-ICNN~\cite{MKLGSB21}.

From the optimality condition, the following relationship holds:
\begin{align}
% Z = \prox_{\eta\tH}(Y) && \iff&& Z - Y = -u, &&& u \in\partial \tH(Z) \label{proxgradL2}\\
Z = \prox_{\eta\Phi}(Y) && \iff&& Z - Y = -\eta u, &&& u \in\partial \Phi (Z). \label{proxgradgeo}
\end{align}
This connects two expressions of the proposed method, a feasible algorithm using the proximal map and an iterative scheme suitable for convergence analysis using the subgradient.

\subsection{Examples}\label{sec:ex}

\subsubsection{Examples of $\FF$}
We give examples of $\FF$ satisfying~\cref{asm:F}, including functionals that are geodesically convex on the constraint but not on the entire Wasserstein space.

Typically, $\FF$ is a potential energy $\EE$, an interaction energy $\WW$, or their sum, where
\begin{equation}
    \EE(\mu) \coloneqq \int V(x) \dd \mu(x), \qquad \WW(\mu) \coloneqq \frac12\iint W(x-y)\dd\mu(x)\dd\mu(y),\label{EW}
\end{equation}
for some $V,W \in C^2(\RR^d)$ such that $W(x) = W(-x)$ for all $x \in \RR^d$.
Their Wasserstein gradients are given by
\[
    \nablaw \EE (\mu) = \nabla V, \qquad \nablaw \WW(\mu) = \int\nabla W(\cdot - y)\dd \mu(y).
\]
If $V$ is $m$-convex, $\EE$ is $m$-geodesically convex. 
However, even if $W$ is $m$-convex with $m > 0$, $\WW$ cannot be $m$-geodesically convex and is only geodesically convex (i.e., $0$-geodesically convex).

\bmhead{Example of $\FF$ which is convex on $\CC$ but nonconvex outside of $\CC$}
Let $d=2$ and for $\mu \in \PP_2(\RR^2)$, define 
$
    (\MM_1(\mu), \MM_2(\mu))^\top \coloneqq \int x \dd \mu(x).
$
We consider 
\begin{equation}
    \FF(\mu)\coloneqq\frac12\int \norm{x}^2\dd\mu(x) + 2\cos\paren*{\mathcal M_1(\mu)},\qquad \CC\coloneqq\setI*{\mu\in\PP_2(\RR^d)}{2\MM_1(\mu) + 1= \MM_2(\mu)}.
\end{equation}
For $X \in L^2(\nu)$, we introduce the notation $\bar X = (\bar X_1, \bar X_2)^\top \coloneqq \int X(z)\dd\nu(z)$. Then, the lifted functional and constraint are
\[
    \tF(X)=\frac12\normL{X}{\nu}^2+2\cos(\bar X_1), \qquad \CCnu = \setI{X\in L^2(\nu)}{2\bar X_1 + 1 = \bar X_2},
\]
implying that $\CCnu$ is closed and convex.
The functional $\FF$ is finite and Wasserstein differentiable on $\PP_2(\RR^d)$, with
\[
    \nablaw\FF(\mu)(x)=x-2\sin(\MM_1(\mu))e_1,\qquad \nabla\tF(X)=X-2\sin(\bar X_1)e_1,
\]
where $e_1 = (1,0)^\top$.

We show that $\FF$ is geodesically convex on $\CC$.
For any $X,Y \in \CCnu$, $T \coloneqq X - Y$ satisfies $2\bar T_1 = \bar T_2$.
Since $\normL{T}{\nu}^2 = \normL{T - \bar T}{\nu}^2 + \norm{\bar T}^2 \ge \norm{\bar T}^2 = 5\bar T_1^2$, we have 
\begin{align}
    \MoveEqLeft
    \inprL{\nabla \tF (X) - \nabla \tF(Y)}{X-Y}{\nu} = \inprL{T - 2(\sin(\bar{X}_1) - \sin(\bar Y_1))e_1}{T}{\nu} \label{ex1m}\\
    &= \normL{T}{\nu}^2 - 2(\sin(\bar{X}_1) - \sin(\bar Y_1))\bar T_1
    \ge \normL{T}{\nu}^2 - 2\bar T_1^2
    \ge \frac35\normL{T}{\nu}^2.
\end{align}
Thus, $\tF$ is $3/5$-convex on $\CCnu$, and by~\cref{prop:convex}, $\FF$ is $3/5$-geodesically convex on $\CC$.

We show that $\FF$ is not geodesically convex on any open neighborhood of $\CC$.
For $s > 0$, consider $\mu_0= \delta_{(0,1)} \in \CC$ and $\mu_s = \delta_{(s,1)} \not\in \CC$. For any open $\mathcal V \supseteq \CC$, by taking sufficiently small $s$, we have $\mu_s \in \mathcal V$ since $W_2(\mu_0,\mu_s) = s$.
For such $\mu_s$, we have
$
    \FF(\mu_s) = \frac12\paren*{s^2+1} + 2 \cos(s),
$
which is nonconvex around $s = 0$.
Since the Wasserstein geodesic from $\mu_0$ to $\mu_s$ is $\mu_{ts}$ for $t \in [0,1]$, this implies that $\FF$ is not geodesically convex on $\mathcal V$.

By a calculation similar to that in~\eqref{ex1m}, we can show that $\FF$ satisfies \cref{asm:F}.\ref{asmLsmooth} with $L = 7/5$.
This is smaller than the constant $L = 3$ in the corresponding unconstrained case.

\bmhead{$m$-geodesic convexity of $\WW$ under a mean constraint}
We consider the interaction energy $\WW$ with $m$-convex $W$ for $m > 0$, defined by~\eqref{EW}. Then, $\WW$ is translation invariant:
\[
    \WW((\id+a)_\#\rho)=\WW(\rho)\qquad\text{for any }a\in\RR^d.
\]
This implies that even though $W$ is $m$-convex, $\WW$ cannot be $m$-geodesically convex.

However, $\WW$ is $m$-geodesically convex on the constraint
$\setI{\mu \in \PP_2(\RR^d)}{\int x \dd \mu(x) = b}$ for $b \in \RR^d$.

% \bmhead{An interaction kernel that is convex only on the relevant difference set.}
% For a concrete example in which global geodesic convexity fails, take $d=1$ and
% \[
%     W(u)=\frac{u^2}{2}+2(1-\cos u),\qquad S=\left[-\frac{\pi}{4},\frac{\pi}{4}\right],\qquad b=0.
% \]
% The kernel is even, its gradient is globally $3$-Lipschitz, and
% \[
%     1\leq W''(u)=1+2\cos u\leq3\qquad\text{for }u\in S-S=\left[-\frac{\pi}{2},\frac{\pi}{2}\right].
% \]
% Hence, on the constraint $\CC_b$,~\cref{asm:F} holds with $m=1$ and $L=3$.

% However, $\FF$ is not geodesically convex on $\PP_2(\RR)$, even if the mean is fixed at zero.
% Indeed, for $a>0$, define
% \[
%     \rho_a\coloneqq\frac12\delta_{-a}+\frac12\delta_a.
% \]
% The curve $a\mapsto\rho_a$ is a constant-speed Wasserstein geodesic on every interval contained in $(0,\infty)$, and
% \[
%     \FF(\rho_a)=\frac14W(2a),\qquad \left.\frac{d^2}{da^2}\FF(\rho_a)\right|_{a=\pi/2}=W''(\pi)=-1.
% \]
% Thus the support constraint excludes regions where the interaction energy has negative curvature, while the prescribed mean removes its translation invariance.

\subsubsection{Example of $\HH$}
We give an example of a regularizer $\HHr$ satisfying~\cref{asm:H}, which can be combined with the constraints discussed above.

\bmhead{Entropy regularization}
For $\beta>0$, consider
\begin{equation}
    \HHr(\mu)=\beta\ent(\mu),\qquad \ent(\mu)\coloneqq\begin{cases}\displaystyle\int p(x)\log p(x)\dd x,&\mu=p\dd x,\\+\infty,&\mu\not\ll\dd x,\end{cases}\label{eq:H-entropy}
\end{equation}
with the convention $0\log0=0$.
The entropy is $W_2$-lower semicontinuous and convex along generalized geodesics, but $\tH$ is not convex in general.
The condition $D(\HH)\neq\varnothing$ requires that $\CC$ contain a measure with finite entropy.

The proximal map is
\begin{equation}
    \prox_{\eta\Phi}(Y)=\argmin_{X\in\CCnu\cap\KK_\nu}\paren*{\beta\ent(X_\#\nu)+\frac{1}{2\eta}\normL{X-Y}{\nu}^2}.\label{entprox}
\end{equation}
For the unconstrained problem $\CC = \PP_2(\RR^d)$, this can be approximated by ICNN.

Although computing the proximal map is generally difficult for $\CC \neq \PP_2(\RR^d)$, this is tractable in some cases.
For example, the mean constraint $\setI{\mu \in \PP_2(\RR^d)}{\int x \dd \mu(x) \in B}$ for a nonempty closed convex set $B \subseteq \RR^d$ can be handled by solving the unconstrained problem and projecting its mean onto $B$. In fact, it holds that
\begin{equation}
    \prox_{\eta\Phi}(Y)=\argmin_{X\in \KK_\nu}\paren*{\beta\ent(X_\#\nu)+\frac{1}{2\eta}\normL{X-Y}{\nu}^2} + P_B\paren*{\int Y\dd\nu}-\int Y\dd\nu,\label{meanentprox}
\end{equation}
where $P_B$ denotes the projection onto $B$.
This is shown in \cref{app:meanconst}.

% For another example, consider the variance constraint~\eqref{variance} and for convex $U \in C^1(\RR)$, set
% \begin{equation}
%     \FF(\mu) = \frac m2 \int \norm{x}^2 \dd \mu(x) + U\paren*{\int x \dd \mu(x)}, \quad \HH(\mu) = \beta \ent(\mu).
% \end{equation}

\section{Algorithms}
In this section, we describe the proposed algorithms for~\eqref{prob} expressed by updates of transport maps and present particle approximations of these algorithms.
For the algorithmic results, we assume that the reference measure $\nu$ belongs to $D(\HH)$ so that we can use the simple initialization $\mu_0=\nu$.

\subsection{Wasserstein FISTA}
Let $L$ and $m$ in \cref{asm:F} satisfy $0\le m<L$ and let $\Phi$ be defined by~\eqref{phi}. 
Define $\setE{a_k}_{k\ge0}$ by $a_0 = 0$ and
\begin{equation}
    a_{k+1}\coloneqq\frac{1+2La_k+\sqrt{1+4La_k(1+ma_k)}}{2(L-m)}.\label{FISTAa}
\end{equation}
Let $\mu_0 = \nu \in \PP_2^{\mathrm{ac}}(\RR^d) \cap \CC\cap D(\HH)$ and $X_0=V_0 = \id$. The update with respect to transport maps is given by
\begin{subequations}\label{FISTA}
\begin{empheq}[left=\empheqlbrace]{align}
    X_{k+1} - X_k &= \frac{a_{k+1}-a_k}{a_k} (V_{k+1} - X_{k+1}), \label{FISTA1}\\
    V_{k+1} - V_k &= \frac{a_{k+1}-a_k}{1 + ma_k} (m(W_{k} - V_{k+1}) - \nabla \tF(W_k) - G_{k+1}), \label{FISTA2}\\
    &\hspace{100pt} G_{k+1} \in \partial \Phi(V_{k+1}), \label{Gupdate}\\
    W_{k} &= \frac{(a_{k+1} - a_k)(1+ma_k)}{a_{k+1}+ 2ma_ka_{k+1} - ma_k^2} V_k + \frac{a_k(1+ma_{k+1})}{a_{k+1}+ 2ma_ka_{k+1} - ma_k^2}X_k. \label{FISTA3}
\end{empheq}
\end{subequations}
When $a_k = 0$, we interpret~\eqref{FISTA1} as $X_{k+1} = V_{k+1}$.
Then, the output is defined by $\mu_k = {X_k}_\#\nu$.
\begin{remark}\label{tem:typical}
    A typical accelerated Forward-backward method in Euclidean space uses an update that corresponds to $G_{k+1} \in \partial \Phi(X_{k+1})$ in~\eqref{Gupdate}. 
    However, this update does not in general preserve $V_k\in\KK_\nu$, as required for our use of~\eqref{convexG}.
\end{remark}
% \begin{remark}
%     We can rewrite the second line of the update as
%     \begin{align}
%         V_{k+1} &= \argmin_{U \in L^2(\nu)} \paren*{ (a_{k+1} - a_k) (\tH(U) + \iota_{\mathcal{K}_\nu}(U)) + \frac12\normL{V_k - (a_{k+1} - a_k)\nabla \tF(W_k) - U}{\nu}^2 }\\
%         &= \argmin_{U \in \mathcal{K}_\nu} \paren*{ (a_{k+1} - a_k) \tH(U)  + \frac12\normL{V_k - (a_{k+1} - a_k)\nabla \tF(W_k) - U}{\nu}^2}\\
%         &= \argmin_{\phi \text{: convex}} \paren*{ (a_{k+1} - a_k) \HH((\nabla \phi)_\# \nu)  + \frac12\normL{V_k - (a_{k+1} - a_k)\nabla \tF(W_k) - \nabla \phi}{\nu}^2}.
%         % &= P_{\mathcal{K}_\nu} \prox_{(a_{k+1} - a_k) \HH}(V_k - (a_{k+1} - a_k)\nablaw \FF(W_k)).
%     \end{align}
%     This subproblem can be handled using the same ICNN parameterization and optimization techniques as those used in .
% \end{remark}
Using the proximal map~\eqref{geoprox} and applying~\eqref{proxgradgeo}, we obtain the explicit form of Wasserstein FISTA.
For $\setE{a_k}_{k\ge0}$ defined by~\eqref{FISTAa}, setting
\begin{equation}
    \tau_k\coloneqq\frac{(a_{k+1}-a_k)(1+ma_k)}{a_{k+1}+2ma_ka_{k+1}-ma_k^2},\qquad \eta_k\coloneqq\frac{a_{k+1}-a_k}{1+ma_{k+1}},\label{FISTAcoef}
\end{equation}
we can write the updates as follows:
\begin{subequations}
\begin{empheq}[left=\empheqlbrace]{align}
    W_k &= (1-\tau_k)X_k+\tau_kV_k,\label{FISTAex1}\\
    Y_{k+1} &= V_k+m\eta_k(W_k - V_k)-\eta_k\nabla\tF(W_k),\label{FISTAexY}\\
    V_{k+1} &= \prox_{\eta_k \Phi}(Y_{k+1}),\label{FISTAex2}\\
    X_{k+1} &= \frac{a_k}{a_{k+1}}X_k+\frac{a_{k+1}-a_k}{a_{k+1}}V_{k+1}.\label{FISTAex3}
\end{empheq}
\end{subequations}

\bmhead{Particle approximation using ICNN}
To implement FISTA numerically, we consider a particle approximation of the algorithm.
Let $x_0^1, \ldots, x_0^N \in \RR^d$ be particles sampled from $\mu_0=\nu$, and approximate $\mu_0$ by $\mu_0^N \coloneqq \frac1N\sum_{i=1}^N \delta_{x_0^i}$.
Define the particles at iteration $k$ by
\[
    x_k^i = X_k(x_0^i), \quad  v_k^i = V_k(x_0^i), \quad w_k^i = W_k(x_0^i), \quad y_k^i = Y_k(x_0^i)
\]
and approximate the output $\mu_k = {X_k}_\#\mu_0$ by $\mu_k^N \coloneqq {X_k}_\#\mu_0^N = \frac1N\sum_{i=1}^N \delta_{x_k^i}$.
Hereafter, we derive an update scheme for $x_k^i, v_k^i, w_k^i$, and $y_{k+1}^i$.

Since $\nabla\tF(W_k) = \nablaw \FF({W_k}_\#\mu_0)\circ W_k$, we approximate $\nabla\tF(W_k) \in L^2(\mu_0)$ by
\begin{align}
    \MoveEqLeft
    \nabla\tF(W_k)(x_0^i) = \nablaw \FF({W_k}_\#\mu_0)(w_k^i)\\
    &\approx \nablaw \FF({W_k}_\#\mu_0^N)(w_k^i) = \nablaw \FF\paren*{\frac1N \sum_{j=1}^N \delta_{w_k^j}}(w_k^i).
\end{align}
Therefore, since
\begin{align}
    {W_k}_\# \mu_0 &= \paren*{(1-\tau_k)X_k+\tau_kV_k}_\#\mu_0,\\
    {Y_{k+1}}_\# \mu_0 &= \paren*{V_k+m\eta_k(W_k - V_k)-\eta_k\nabla\tF(W_k)}_\# \mu_0,\\
    {X_{k+1}}_\#\mu_0 &= \paren*{\frac{a_k}{a_{k+1}}X_k+\frac{a_{k+1}-a_k}{a_{k+1}}V_{k+1}}_\#\mu_0,
\end{align}
from~\eqref{FISTAex1},~\eqref{FISTAexY} and~\eqref{FISTAex3},
we have particle updates
\begin{align}
    w_k^i &= (1-\tau_k)x_k^i+\tau_kv_k^i,\\
    y_{k+1}^i &= v_k^i+m\eta_k(w_k^i - v_k^i)-\eta_k\nablaw \FF\paren*{\frac1N \sum_{j=1}^N \delta_{w_k^j}}(w_k^i),\\
    x_{k+1}^i &= \frac{a_k}{a_{k+1}}x_k^i+\frac{a_{k+1}-a_k}{a_{k+1}}v_{k+1}^i.
\end{align}
Then, we derive a particle approximation for~\eqref{FISTAex2} using the expression~\eqref{proxB}.
Since
\[
    \normL{\nabla \phi-Y_{k+1}}{\mu_0}^2 = \int \norm{\nabla \phi (x)-Y_{k+1}(x)}^2 \dd \mu_0(x) \approx \frac1N \sum_{i=1}^N \norm{\nabla \phi (x_0^i)-y_{k+1}^i}^2,
\]
we obtain the approximation
\begin{align}
    \prox_{\eta_k\Phi}(Y_{k+1}) &\approx \nabla \hat\phi_{Y_{k+1}},\\
    \hat\phi_{Y_{k+1}} &\in \argmin_{\substack{\phi \colon \text{convex}\\ \text{s.t. } \nabla \phi \in \CCnu}}\paren*{\tHr^N(\nabla \phi)+\frac{1}{2\eta_k N} \sum_{i=1}^N \norm{\nabla \phi (x_0^i)-y_{k+1}^i}^2}, \label{proxapprox}
\end{align}
where $\tHr^N(\nabla \phi)$ is an approximation of $\tHr(\nabla \phi)$ constructed from the reference samples $x_0^1,\ldots,x_0^N$.
Such $\hat\phi_{Y_{k+1}}$ can be approximately constructed using ICNNs, as discussed in \cref{sec:prox}.
With this approximation, we approximate~\eqref{FISTAex2} by setting
\[
    v_{k+1}^i = \nabla \hat\phi_{Y_{k+1}}(x_0^i).
\]
Combining these updates, we obtain~\cref{alg:FISTA}.

\begin{algorithm}[htbp]
    \caption{Particle implementation of FISTA}
    \label{alg:FISTA}
    \begin{algorithmic}[1]
        \State \textbf{Input:} $L$, $m$, Initial particles $\mathbf x_0 = (x_0^1,\ldots,x_0^N) \in (\RR^d)^N$, $K$
        \State Initialize $\mathbf v_0 = (v_0^1,\ldots,v_0^N) \gets \mathbf x_0$, $\mathbf w_0 = (w_0^1,\ldots,w_0^N) \gets \mathbf x_0$.
        \For{$k=0,\ldots,K-1$}
            \State $\mathbf w_k \gets(1-\tau_k)\mathbf x_k+\tau_k \mathbf v_k$
            \Comment{$\tau_k$ is defined by~\eqref{FISTAcoef}}
            \State $y_{k+1}^i= v_k^i+m\eta_k(w_k^i - v_k^i)-\eta_k\nablaw \FF\paren*{\frac{1}{N}\sum_{j=1}^N \delta_{w_k^j}}(w_k^i)$ \quad for $i = 1,\ldots,N$\\
            \Comment{$\eta_k$ is defined by~\eqref{FISTAcoef}}
            \State $v_{k+1}^i\gets \mathsf P_{k}(x_0^i)$ \quad for $i = 1,\ldots,N$\\
            \Comment{$\mathsf P_{k}$ approximates $\prox_{\eta_k\Phi}(Y_{k+1})$, e.g., via~\eqref{proxapprox}}
            \State $\mathbf x_{k+1}\gets\displaystyle\frac{a_k}{a_{k+1}}\mathbf x_k+\frac{a_{k+1}-a_k}{a_{k+1}}\mathbf v_{k+1}$.
        \EndFor
        \State \textbf{Output:} $\displaystyle\mu_K^N=\frac1N\sum_{i=1}^N\delta_{x_K^i}$.
    \end{algorithmic}
\end{algorithm}

\subsection{Wasserstein SR2FISTA}
We consider another accelerated proximal gradient method with a better asymptotic linear convergence factor when $m>0$.
Define $\setE{a_k}_{k\ge0}$ by $a_0 = 0$ and 
\begin{equation}
a_{k+1}\coloneqq\frac{1+La_k+\sqrt{(1+ma_k)(1+(2L-m)a_k)}}{L-m} \label{SR2FISTAa}
\end{equation}
Let $\mu_0 = \nu \in \PP_2^{\mathrm{ac}}(\RR^d) \cap \CC \cap D(\HH)$ and $X_0=V_0 = \id$. The update with respect to transport maps is given by
\begin{subequations}\label{SR2FISTA}
\begin{empheq}[left=\empheqlbrace]{align}
        X_{k+1} - X_k &= \frac{a_{k+1}-a_k}{a_k} (V_{k+1} - X_{k+1}),\label{TWNAG1}\\
        V_{k+1} - V_k &= \frac{a_{k+1}-a_k}{2(1 + ma_k)} (m(W_{k} - V_{k+1}) - \nabla \tF(W_k) - G_{k+1}), \label{TWNAG2}\\
        &\hspace{100pt} G_{k+1} \in \partial \Phi(V_{k+1}),\\
        W_{k} &= \frac{a_{k+1}-a_k}{a_{k+1}} V_k + \frac{a_k}{a_{k+1}}X_k, \label{TWNAG3}
\end{empheq}
\end{subequations}
When $a_k = 0$, we interpret~\eqref{TWNAG1} as $X_{k+1} = V_{k+1}$.
Then, $\mu_k = {X_k}_\#\nu$.

For $\setE{a_k}_{k\ge0}$ defined by~\eqref{SR2FISTAa}, setting
\begin{equation}
    \tau_k\coloneqq\frac{a_{k+1}-a_k}{a_{k+1}},\qquad \eta_k\coloneqq\frac{a_{k+1}-a_k}{2+m(a_k+a_{k+1})},\label{SR2FISTAcoef}
\end{equation}
we can write~\eqref{SR2FISTA} in an explicit form as
\begin{subequations}\label{SR2FISTAex}
\begin{empheq}[left=\empheqlbrace]{align}
    W_k &= (1-\tau_k)X_k+\tau_kV_k,\label{SR2FISTAex1}\\
    Y_{k+1} &= V_k+ m\eta_k(W_k- V_k)-\eta_k\nabla\tF(W_k),\\
    V_{k+1} &= \prox_{\eta_k\Phi}(Y_{k+1}),\label{SR2FISTAex2}\\
    X_{k+1} &= (1-\tau_k)X_k+\tau_kV_{k+1}.\label{SR2FISTAex3}
\end{empheq}
\end{subequations}

Using the same approximation procedure, we obtain a particle implementation of SR2FISTA, as presented in \cref{alg:SR2FISTA}.
\begin{algorithm}[htbp]
    \caption{Particle implementation of SR2FISTA}
    \label{alg:SR2FISTA}
    \begin{algorithmic}[1]
        \State \textbf{Input:} $L$, $m$, Initial particles $\mathbf x_0 = (x_0^1,\ldots,x_0^N) \in (\RR^d)^N$, $K$
        \State Initialize $\mathbf v_0 = (v_0^1,\ldots,v_0^N) \gets \mathbf x_0$, $\mathbf w_0 = (w_0^1,\ldots,w_0^N) \gets \mathbf x_0$.
        \For{$k=0,\ldots,K-1$}
            \State $\mathbf w_k \gets(1-\tau_k)\mathbf x_k+\tau_k \mathbf v_k$
            \Comment{$\tau_k$ is defined by~\eqref{SR2FISTAcoef}}
            \State $y_{k+1}^i= v_k^i+m\eta_k(w_k^i - v_k^i)-\eta_k\nablaw \FF\paren*{\frac{1}{N}\sum_{j=1}^N \delta_{w_k^j}}(w_k^i)$ \quad for $i = 1,\ldots,N$\\
            \Comment{$\eta_k$ is defined by~\eqref{SR2FISTAcoef}}
            \State $v_{k+1}^i\gets \mathsf P_{k}(x_0^i)$ \quad for $i = 1,\ldots,N$\\
            \Comment{$\mathsf P_{k}$ approximates $\prox_{\eta_k\Phi}(Y_{k+1})$, e.g., via~\eqref{proxapprox}}
            \State $\mathbf x_{k+1}\gets\displaystyle (1 - \tau_k)\mathbf x_k+\tau_k\mathbf v_{k+1}$.
        \EndFor
        \State \textbf{Output:} $\displaystyle\mu_K^N=\frac1N\sum_{i=1}^N\delta_{x_K^i}$.
    \end{algorithmic}
\end{algorithm}

\section{Convergence analyses}\label{sec:convergence}
In this section, we establish convergence rates for the exact transport-map iterations~\eqref{FISTA} and~\eqref{SR2FISTA}.
We begin by showing lemmas that are commonly used in the analyses of both algorithms.

\begin{lemma}\label{lem:opt}
    For Wasserstein FISTA~\eqref{FISTA} and SR2FISTA~\eqref{SR2FISTA}, we have $X_{k}, V_k,W_k \in \KK_\nu\cap\CCnu$ for all $k \ge 0$.
\end{lemma}
\begin{proof}
    By the proximal update, $V_{k+1}\in D(\Phi)\subseteq\KK_\nu\cap\CCnu$.
    Since $X_0=V_0=\id\in\KK_\nu\cap\CCnu$, we proceed by induction and suppose $X_{k} \in \KK_\nu\cap\CCnu$. Since $\KK_\nu\cap\CCnu$ is convex, it follows from~\eqref{FISTAex1} and \eqref{FISTAex3} for FISTA, and from \eqref{SR2FISTAex1} and \eqref{SR2FISTAex3} for SR2FISTA that $X_{k+1},W_{k} \in \KK_\nu\cap\CCnu$.
\end{proof}

\begin{lemma}\label{lem:norm}
    For $X_{k+1},X_k,V_{k+1}$ satisfying \eqref{FISTA1}, or equivalently~\eqref{TWNAG1}, it holds that for any $Z \in L^2(\nu)$,
    \begin{align}
        \MoveEqLeft
        a_k\normL{X_k - Z}{\nu}^2 + (a_{k+1} - a_k)\normL{V_{k+1} - Z}{\nu}^2 - a_{k+1}\normL{X_{k+1}-Z}{\nu}^2\\
        &= \frac{a_k(a_{k+1} - a_k)}{a_{k+1}}\normL{X_k - V_{k+1}}{\nu}^2.
    \end{align}
\end{lemma}
\begin{proof}
    We use the equality
    \begin{equation}
        \normL{\alpha S + (1-\alpha) T}{\nu}^2 = \alpha \normL{S}{\nu}^2 + (1-\alpha) \normL{T}{\nu}^2 - \alpha(1-\alpha)\normL{S-T}{\nu}^2 \label{alpha}
    \end{equation}
    for any $S,T \in L^2(\nu)$ and $\alpha \in \RR$. This, together with~\eqref{FISTAex3}, or equivalently~\eqref{SR2FISTAex3}, yields for any $Z \in L^2(\nu)$, 
    \begin{align}
    \MoveEqLeft
        \normL{X_{k+1}-Z}{\nu}^2\\
        &=\normL*{\frac{a_k}{a_{k+1}}(X_k - Z) + \frac{a_{k+1} - a_k}{a_{k+1}}(V_{k+1} - Z)}{\nu}^2\\
        &=\frac{a_k}{a_{k+1}}\normL{X_k - Z}{\nu}^2 + \frac{a_{k+1} - a_k}{a_{k+1}}\normL{V_{k+1} - Z}{\nu}^2 - \frac{a_k(a_{k+1} - a_k)}{a_{k+1}^2}\normL{X_k - V_{k+1}}{\nu}^2.
    \end{align}
\end{proof}
In the following lemma, we show four inequalities. 
For Euclidean FISTA, the corresponding inequalities in $\RR^d$ are the sufficient ingredients for the convergence analysis.
By showing the inequalities in $L^2(\nu)$, we can carry out essentially the same convergence analysis.
\begin{lemma}\label{lem:ineq}
    Assume \cref{asm:C}, \cref{asm:F} and \cref{asm:H}.
    For FISTA~\eqref{FISTA} and SR2FISTA~\eqref{SR2FISTA}, the following inequalities hold.
    \begin{align}
        &\tF(X_{k+1}) - \tF(X^\star) \le \inprL{\nabla \tF(W_k)}{X_{k+1} - X^\star}{\nu} + \frac L2\normL{X_{k+1} - W_k}{\nu}^2 - \frac m2 \normL{X^\star - W_k}{\nu}^2, \label{Fineq1}\\
        &\tF(X_{k+1}) - \tF(X_k) \le \inprL{\nabla \tF(W_k)}{X_{k+1} - X_k}{\nu} + \frac L2\normL{X_{k+1} - W_k}{\nu}^2 - \frac m2 \normL{X_{k} - W_k}{\nu}^2, \label{Fineq2}\\
        &a_{k+1}\tH(X_{k+1}) - a_k \tH(X_k) \le (a_{k+1} - a_k)\tH(V_{k+1}), \label{jensen}\\
        &\tH(V_{k+1}) - \tH(X^\star) \le \inprL{G_{k+1}}{V_{k+1} - X^\star}{\nu}, \label{Hineq}
    \end{align}
    where $X^\star \in \KK_\nu \cap D(\tH)$.
    % Moreover, additionally assume \cref{asm:H2}. Then, the above inequalities hold for FISTA~\eqref{FISTA} and SR2FISTA~\eqref{SR2FISTA} with~\eqref{GL2}.
\end{lemma}
\begin{proof}
    For $\tF$, it holds that for any $S,T,W \in \CCnu$,
    \begin{equation}
        \tF(S) - \tF(T) \le \inprL{\nabla\tF(W)}{S-T}{\nu} + \frac L2 \normL{S-W}{\nu}^2 - \frac m2 \normL{T-W}{\nu}^2, \label{3pt}
    \end{equation}
    which follows by summing the two inequalities from~\eqref{Lsmooth} and~\eqref{convexF}:
    \begin{align}
        \tF(S) - \tF(W) &\le \inprL{\nabla\tF(W)}{S-W}{\nu} + \frac L2 \normL{S-W}{\nu}^2,\\
        \tF(W) - \tF(T) &\le \inprL{\nabla\tF(W)}{W-T}{\nu} - \frac m2 \normL{T-W}{\nu}^2,
    \end{align}
    for any $S,T,W \in \CCnu$.
    By substituting appropriate variables into~\eqref{3pt}, we obtain~\eqref{Fineq1} and~\eqref{Fineq2}.

    Since $\HH$ is convex along generalized geodesics, and since we have~\eqref{FISTAex3} for FISTA and ~\eqref{SR2FISTAex3} for SR2FISTA, together with $X_{k} \in \mathcal{K}_\nu$ and $V_{k+1} \in \mathcal{K}_\nu$,
    it follows that
    \[
        \tH (X_{k+1}) \le \frac{a_k}{a_{k+1}}\tH (X_k) + \frac{a_{k+1} - a_k}{a_{k+1}}\tH(V_{k+1}),
    \]
    which yields~\eqref{jensen}. \eqref{Hineq} follows from~\eqref{convexG}.
\end{proof}

\subsection{Wasserstein FISTA}
\begin{lemma}\label{lem:lyapu1}
    Suppose that $\FF$, $\HH$ and $\CC$ satisfy \cref{asm:C}, \cref{asm:F}, and \cref{asm:H}.
    For $\mu_k$ and $V_k$ of Wasserstein FISTA~\eqref{FISTA}, define the Lyapunov function by
    \begin{align}
    \mathcal{E}_k &= a_k\paren*{\GG(\mu_k) - \GG(\mu^\star)} + \frac{1+ma_k}{2}\normL{V_k - X^\star}{\nu}^2\\
    &= a_k\paren*{\tF(X_k) + \tH(X_k) - \tF(X^\star) - \tH(X^\star)} + \frac{1+ma_k}{2}\normL{V_k - X^\star}{\nu}^2,
    \end{align}
    where $X^\star \in \KK_\nu$ denotes the optimal transport map from $\nu$ to $\mu^\star$.
    Then, $\mathcal{E}_k$ is nonincreasing.
\end{lemma}
From \cref{lem:ineq}, the proof can be carried out in essentially the same manner as the Euclidean case~\cite[Theorem~4.22]{ADA21}. Here, we provide another expression of the proof, which makes the algebraic manipulations explicit.
\begin{proof}
    We have
    \begin{align}
        \MoveEqLeft
        \mathcal{E}_{k+1} - \mathcal{E}_k\\
        &\overset{(\circ)}{\le} (a_{k+1} - a_k)\paren*{\tF(X_{k+1}) - \tF(X^\star)} + a_k\paren*{\tF(X_{k+1}) - \tF(X_k)}\\
        &\qquad + (a_{k+1} - a_k) (\tH(V_{k+1}) - \tH(X^\star))\\
        &\qquad + \frac m2 (a_{k+1} - a_k)\normL{V_{k+1}-X^\star}{\nu}^2 + \frac{1+ma_k}{2}\paren*{\normL{V_{k+1}-X^\star}{\nu}^2 - \normL{V_{k}-X^\star}{\nu}^2}\\
        &\overset{(\bullet)}{\le} (a_{k+1} - a_k) \paren*{\inprL{\nabla \tF(W_k)}{X_{k+1} - X^\star}{\nu} + \frac L2\normL{X_{k+1} - W_k}{\nu}^2 - \frac m2 \normL{X^\star - W_k}{\nu}^2}\\
        &\qquad + a_k\paren*{\inprL{\nabla \tF(W_k)}{X_{k+1} - X_k}{\nu} + \frac L2\normL{X_{k+1} - W_k}{\nu}^2 - \frac m2 \normL{X_{k} - W_k}{\nu}^2}\\
        &\qquad + (a_{k+1} - a_k)\inprL{G_{k+1}}{V_{k+1} - X^\star}{\nu}\\
        &\qquad + \frac m2 (a_{k+1} - a_k)\normL{V_{k+1}-X^\star}{\nu}^2 \\
        &\qquad + \frac{1+ma_k}{2}\paren*{2\inprL{V_{k+1} - V_k}{V_{k+1} - X^\star}{\nu} -\normL{V_{k+1} - V_k}{\nu}^2}\\
        &= (a_{k+1} - a_k) \inprL{\nabla \tF(W_k)}{X_{k+1} - X^\star}{\nu} + \inprL*{\nabla \tF(W_k)}{a_k(X_{k+1} - X_{k})}{\nu} \\
        &\qquad + \inprL{(a_{k+1} - a_k)G_{k+1} + (1 + ma_k)(V_{k+1} - V_k)}{V_{k+1} - X^\star}{\nu} + Q_k,
    \end{align}
    where~\eqref{jensen} is applied to $(\circ)$, \eqref{Fineq1}, \eqref{Fineq2} and \eqref{Hineq} are applied to $(\bullet)$, and we set
    \begin{align}
    \MoveEqLeft Q_k \coloneqq \frac L2 a_{k+1}\normL{X_{k+1} - W_k}{\nu}^2  - \frac m2 (a_{k+1} - a_k)\normL{X^\star - W_k}{\nu}^2 - \frac m2 a_k\normL{X_{k} - W_k}{\nu}^2\\
    &\quad + \frac m2 (a_{k+1} - a_k)\normL{V_{k+1}-X^\star}{\nu}^2 - \frac{1+ma_k}{2}\normL{V_{k+1} - V_k}{\nu}^2.
    \end{align}
    By substituting \eqref{FISTA1} and \eqref{FISTA2} into $X_{k+1} - X_k$ and $V_{k+1} - V_k$, we have
    \begin{align}
    \MoveEqLeft
        \mathcal{E}_{k+1} - \mathcal{E}_k\\
        &\le (a_{k+1} - a_k) \inprL{\nabla \tF(W_k)}{X_{k+1} - X^\star}{\nu} + \inprL*{\nabla \tF(W_k)}{(a_{k+1} - a_k)(V_{k+1} - X_{k+1})}{\nu}\\
        &\qquad + \inprL*{(a_{k+1} - a_k)(m(W_{k} - V_{k+1}) - \nabla \tF(W_k))}{V_{k+1} - X^\star}{\nu} + Q_k\\
        &= m(a_{k+1} - a_k)\inprL*{W_{k} - V_{k+1}}{V_{k+1} - X^\star}{\nu} + Q_k.
    \end{align}
    Since
    \[
        \inprL*{W_{k} - V_{k+1}}{V_{k+1} - X^\star}{\nu} = - \frac12 \normL{W_k - V_{k+1}}{\nu}^2 - \frac12 \normL{V_{k+1} - X^\star}{\nu}^2 + \frac12 \normL{W_k - X^\star}{\nu}^2,
    \]
    whose last two terms are common to $Q_k$, we have
    \begin{align}
        \mathcal{E}_{k+1} - \mathcal{E}_k
        &\le \frac L2 a_{k+1}\normL{X_{k+1} - W_k}{\nu}^2 - \frac m2 a_k\normL{X_{k} - W_k}{\nu}^2\\
        &\quad  - \frac{1+ma_k}{2}\normL{V_{k+1} - V_k}{\nu}^2 - \frac m2 (a_{k+1} - a_k)\normL{W_k - V_{k+1}}{\nu}^2.
    \end{align}
    By substituting $W_k$ into $Z$ of \cref{lem:norm}, we have 
    \begin{align}
        \MoveEqLeft
        - a_k\normL{X_k - W_k}{\nu}^2 - (a_{k+1} - a_k)\normL{V_{k+1} - W_k}{\nu}^2 \\
        &= - a_{k+1}\normL{X_{k+1}- W_k}{\nu}^2 - \frac{a_k(a_{k+1} - a_k)}{a_{k+1}}\normL{X_k - V_{k+1}}{\nu}^2.
    \end{align}
    Therefore, combining the above two expressions, we have
    \begin{align}
        \mathcal{E}_{k+1} - \mathcal{E}_k
        &\le \paren*{\frac L2 - \frac m2} a_{k+1}\normL{X_{k+1} - W_k}{\nu}^2\\
        &\qquad - \frac{1+ma_k}{2}\normL{V_{k+1} - V_k}{\nu}^2 - \frac{ma_k(a_{k+1} - a_k)}{2a_{k+1}}\normL{V_{k+1} - X_k}{\nu}^2.
    \end{align}
    Here, letting
    \[
        \alpha_1 \coloneqq \frac{1+ma_k}{2} + \frac{ma_k(a_{k+1} - a_k)}{2a_{k+1}},\quad  \alpha_V \coloneqq \frac{1}{\alpha_1}\frac{1+ma_k}{2}, \quad \alpha_X \coloneqq \frac{1}{\alpha_1}\frac{ma_k(a_{k+1} - a_k)}{2a_{k+1}},
    \]
    and applying~\eqref{alpha} with $S = V_{k+1} - V_k$ and $T = V_{k+1} - X_k$,
    we have
    \[
        \alpha_V \normL{V_{k+1} - V_k}{\nu}^2 + \alpha_X \normL{V_{k+1} - X_k}{\nu}^2 =  \normL{V_{k+1} - \alpha_V V_k - \alpha_X X_k}{\nu}^2 + \alpha_V\alpha_X\normL{V_k - X_k}{\nu}^2.
    \]
    Thus, discarding the term $\alpha_1\alpha_V\alpha_X\normL{V_k - X_k}{\nu}^2$, we have 
    \[
        \mathcal{E}_{k+1} - \mathcal{E}_k \le \paren*{\frac L2 - \frac m2} a_{k+1}\normL{X_{k+1} - W_k}{\nu}^2  - \alpha_1\normL{V_{k+1} - \alpha_V V_k - \alpha_X X_k}{\nu}^2.
    \]
    From~\eqref{FISTA1}, we have 
    \[
        V_{k+1} = \frac{a_{k+1}}{a_{k+1}-a_k}X_{k+1} - \frac{a_{k}}{a_{k+1}-a_k}X_k,
    \]
    and from~\eqref{FISTA3}, we have
    \[
        W_k = \frac{a_k}{a_{k+1}}X_k + \frac{a_{k+1} - a_k}{a_{k+1}}\paren*{\alpha_V V_k + \alpha_X X_k}.
    \]
    Therefore,
    \begin{align}
    \mathcal{E}_{k+1} - \mathcal{E}_k &\le
        \paren*{\frac L2 - \frac m2}a_{k+1} \normL{X_{k+1} - W_k}{\nu}^2\\
        &\qquad - \alpha_1\normL*{\frac{a_{k+1}}{a_{k+1}-a_k}X_{k+1} - \frac{a_{k}}{a_{k+1}-a_k}X_k - \alpha_V V_k - \alpha_X X_k}{\nu}^2\\
        &= \paren*{\paren*{\frac L2 - \frac m2}a_{k+1} - \alpha_1 \paren*{\frac{a_{k+1}}{a_{k+1}-a_k}}^2} \normL{X_{k+1} - W_k}{\nu}^2,
    \end{align}
    which is nonpositive if
    \[
        L(a_{k+1} - a_k)^2 - a_{k+1} - ma_{k+1}^2 \le 0.
    \]
    Since \eqref{FISTAa} implies $L(a_{k+1} - a_k)^2 - a_{k+1} - ma_{k+1}^2 = 0$, the claim follows.
\end{proof}

\begin{lemma}\label{lem:FISTAa}
    For $a_k$ defined by~\eqref{FISTAa}, the following bound holds for any $k\ge1$:
    \[
        \frac 1{a_k} \le \min\setE*{\frac{2}{k^2},\paren*{1-\sqrt{\frac mL}}^k} 2L.
    \]
\end{lemma}
The proof is given in~\cite[Corollary~4.13]{ADA21}. For completeness, we give a proof in \cref{app:convergence}.

\begin{theorem}[Convergence rate of Wasserstein FISTA]\label{thm:FISTA}
    Assume \cref{asm:C}, \cref{asm:F}, and \cref{asm:H}. Then, for Wasserstein FISTA~\eqref{FISTA}, the following bound holds for any $k\ge1$:
    \[
        \GG(\mu_k) - \GG(\mu^\star) \le \min\setE*{\frac{2}{k^2},\paren*{1-\sqrt{\frac mL}}^k} LW_2(\mu_0,\mu^\star)^2.
    \]
\end{theorem}
\begin{proof}
    From \cref{lem:lyapu1}, we have for any $k \ge 0$
    \[
        a_k(\GG(\mu_k) - \GG(\mu^\star)) \le \mathcal E_k \le \mathcal E_0 = \frac12 \normL{V_0 - X^\star}{\nu}^2 = \frac 12 \normL{\id - X^\star}{\mu_0}^2 = \frac12 W_2^2(\mu_0,\mu^\star).
    \]
    Therefore, from \cref{lem:FISTAa}, the claim follows.
\end{proof}

\subsection{Wasserstein SR2FISTA}

\begin{lemma}\label{lem:lyapu2}
    Suppose that $\FF$, $\HH$ and $\CC$ satisfy \cref{asm:C}, \cref{asm:F}, and \cref{asm:H}.
    For $\mu_k$ and $V_k$ of Wasserstein SR2FISTA~\eqref{SR2FISTA}, define the Lyapunov function by
    \begin{align}
    \mathcal{E}_k &= a_k\paren*{\GG(\mu_k) - \GG(\mu^\star) - \frac m2 \normL{X_k - X^\star}{\nu}^2} + (1+ma_k)\normL{V_k - X^\star}{\nu}^2\\
    &= a_k\paren*{\tF(X_k) + \tH(X_k) - \tF(X^\star) - \tH(X^\star) - \frac m2 \normL{X_k - X^\star}{\nu}^2} + (1+ma_k)\normL{V_k - X^\star}{\nu}^2, \label{lyapusr2}
    \end{align}
    where $X^\star \in \KK_\nu$ denotes the optimal transport map from $\nu$ to $\mu^\star$.
    Then, $\mathcal{E}_k$ is nonincreasing.
\end{lemma}
\begin{remark}
    This Lyapunov function has essentially the same form as that used for Euclidean SR2FISTA~\cite{U25}, which adopts the typical forward-backward scheme corresponding to $G_{k+1} \in \partial \Phi(X_{k+1})$ as mentioned in \cref{tem:typical}.
    This lemma shows that the Lyapunov function~\eqref{lyapusr2} is still nonincreasing for the variant of SR2FISTA considered here.
\end{remark}
\begin{proof}    
    We use the inequalities shown in~\cref{lem:ineq} as in the FISTA case.
    \begin{align}
        \MoveEqLeft
        \mathcal{E}_{k+1} - \mathcal{E}_k\\
        &\overset{(\circ)}{\le} (a_{k+1} - a_k)\paren*{\tF(X_{k+1}) - \tF(X^\star)} + a_k\paren*{\tF(X_{k+1}) - \tF(X_k)}\\
        &\qquad + (a_{k+1} - a_k) (\tH(V_{k+1}) - \tH(X^\star))\\
        &\qquad - \frac m2 \paren*{a_{k+1}\normL{X_{k+1}-X^\star}{\nu}^2 -a_k\normL{X_{k}-X^\star}{\nu}^2} \\
        &\qquad + m (a_{k+1} - a_k)\normL{V_{k+1}-X^\star}{\nu}^2 + (1+ma_k)\paren*{\normL{V_{k+1}-X^\star}{\nu}^2 - \normL{V_{k}-X^\star}{\nu}^2}\\
        &\overset{(\bullet)}{\le} (a_{k+1} - a_k) \paren*{\inprL{\nabla \tF(W_k)}{X_{k+1} - X^\star}{\nu} + \frac L2\normL{X_{k+1} - W_k}{\nu}^2 - \frac m2 \normL{X^\star - W_k}{\nu}^2}\\
        &\qquad + a_k\paren*{\inprL{\nabla \tF(W_k)}{X_{k+1} - X_k}{\nu} + \frac L2\normL{X_{k+1} - W_k}{\nu}^2 - \frac m2 \normL{X_{k} - W_k}{\nu}^2}\\
        &\qquad + (a_{k+1} - a_k)\inprL{G_{k+1}}{V_{k+1} - X^\star}{\nu}\\
        &\qquad - \frac m2 \paren*{a_{k+1}\normL{X_{k+1}-X^\star}{\nu}^2 -a_k\normL{X_{k}-X^\star}{\nu}^2} + m (a_{k+1} - a_k)\normL{V_{k+1}-X^\star}{\nu}^2 \\
        &\qquad + (1+ma_k)\paren*{2\inprL{V_{k+1} - V_k}{V_{k+1} - X^\star}{\nu} -\normL{V_{k+1} - V_k}{\nu}^2}\\
        &= (a_{k+1} - a_k) \inprL{\nabla \tF(W_k)}{X_{k+1} - X^\star}{\nu} + \inprL*{\nabla \tF(W_k)}{a_k(X_{k+1} - X_{k})}{\nu} \\
        &\qquad + \inprL{(a_{k+1} - a_k)G_{k+1} + 2(1 + ma_k)(V_{k+1} - V_k)}{V_{k+1} - X^\star}{\nu} + Q_k,
    \end{align}
    where~\eqref{jensen} is applied to $(\circ)$, \eqref{Fineq1}, \eqref{Fineq2} and \eqref{Hineq} are applied to $(\bullet)$, and we set
    \begin{align}
    \MoveEqLeft Q_k \coloneqq \frac L2 a_{k+1}\normL{X_{k+1} - W_k}{\nu}^2  - \frac m2 (a_{k+1} - a_k)\normL{X^\star - W_k}{\nu}^2 - \frac m2 a_k\normL{X_{k} - W_k}{\nu}^2\\
    &\quad - \frac m2 \paren*{a_{k+1}\normL{X_{k+1}-X^\star}{\nu}^2 -a_k\normL{X_{k}-X^\star}{\nu}^2} + m (a_{k+1} - a_k)\normL{V_{k+1}-X^\star}{\nu}^2\\
    &\quad  - (1+ma_k)\normL{V_{k+1} - V_k}{\nu}^2.
    \end{align}
    By substituting \eqref{TWNAG1} and \eqref{TWNAG2} into $X_{k+1} - X_k$ and $V_{k+1} - V_k$, we have
    \begin{align}
    \MoveEqLeft
        \mathcal{E}_{k+1} - \mathcal{E}_k\\
        &\le (a_{k+1} - a_k) \inprL{\nabla \tF(W_k)}{X_{k+1} - X^\star}{\nu} + \inprL*{\nabla \tF(W_k)}{(a_{k+1} - a_k)(V_{k+1} - X_{k+1})}{\nu}\\
        &\qquad + \inprL*{(a_{k+1} - a_k)(m(W_{k} - V_{k+1}) - \nabla \tF(W_k))}{V_{k+1} - X^\star}{\nu} + Q_k\\
        &= m(a_{k+1} - a_k)\inprL*{W_{k} - V_{k+1}}{V_{k+1} - X^\star}{\nu} + Q_k.
    \end{align}
    Since
    \[
        \inprL*{W_{k} - V_{k+1}}{V_{k+1} - X^\star}{\nu} = - \frac12 \normL{W_k - V_{k+1}}{\nu}^2 - \frac12 \normL{V_{k+1} - X^\star}{\nu}^2 + \frac12 \normL{W_k - X^\star}{\nu}^2,
    \]
    whose last two terms are common to $Q_k$, we have
    \begin{align}
        \MoveEqLeft
        \mathcal{E}_{k+1} - \mathcal{E}_k\\
        &\le \frac L2 a_{k+1}\normL{X_{k+1} - W_k}{\nu}^2 - \frac m2 a_k\normL{X_{k} - W_k}{\nu}^2\\
        &\quad - \frac m2 \paren*{a_{k+1}\normL{X_{k+1}-X^\star}{\nu}^2 -a_k\normL{X_{k}-X^\star}{\nu}^2} + \frac m2 (a_{k+1} - a_k)\normL{V_{k+1}-X^\star}{\nu}^2\\
        &\quad  - (1+ma_k)\normL{V_{k+1} - V_k}{\nu}^2 - \frac m2 (a_{k+1} - a_k)\normL{W_k - V_{k+1}}{\nu}^2.\label{Emid}
    \end{align}
    By substituting $X^\star$ and $W_k$ into $Z$ of \cref{lem:norm}, we have 
    \begin{align}
        &-a_{k+1} \normL{X_{k+1}-X^\star}{\nu}^2 + a_k\normL{X_k - X^\star}{\nu}^2 + (a_{k+1} - a_k)\normL{V_{k+1} - X^\star}{\nu}^2 \\
        =& -a_{k+1} \normL{X_{k+1}-W_k}{\nu}^2 + a_k\normL{X_k -W_k}{\nu}^2 + (a_{k+1} - a_k)\normL{V_{k+1} -W_k}{\nu}^2.
    \end{align}
    Therefore, applying this equality to~\eqref{Emid} and using~\eqref{TWNAG1} and \eqref{TWNAG3}, we have
    \begin{align}
    \MoveEqLeft
        \mathcal{E}_{k+1} - \mathcal{E}_k \le \paren*{\frac L2 - \frac m2} a_{k+1}\normL{X_{k+1} - W_k}{\nu}^2  - (1+ma_k)\normL{V_{k+1} - V_k}{\nu}^2\\
        &= \paren*{\frac L2 - \frac m2}a_{k+1} \normL{X_{k+1} - W_k}{\nu}^2  - (1+ma_k)\normL*{\frac{a_{k+1}}{a_{k+1}-a_k}X_{k+1} - \frac{a_{k}}{a_{k+1}-a_k}X_k - V_k}{\nu}^2\\
        &= \paren*{\paren*{\frac L2 - \frac m2}a_{k+1} - (1+ma_k) \paren*{\frac{a_{k+1}}{a_{k+1}-a_k}}^2} \normL{X_{k+1} - W_k}{\nu}^2,
    \end{align}
    which is nonpositive if
    \[
        (L-m)(a_{k+1} - a_k)^2 - 2(1+ma_k)a_{k+1} \le 0.
    \]
    Since \eqref{SR2FISTAa} implies that $(L-m)(a_{k+1} - a_k)^2 - 2(1+ma_k)a_{k+1} = 0$, the claim follows.
\end{proof}

\begin{lemma}\label{lem:SR2FISTAa}
    For $a_k$ defined by~\eqref{SR2FISTAa}, the following bound holds for any $k\ge1$:
    \begin{equation}
        \frac 1{a_k} \le \min\setE*{\frac{2L}{k^2},\frac{L-m}{2}r^{-(k-1)}}, \quad \text{where } r \coloneqq \frac{1 + \sqrt{\frac{m}{L}\paren*{2 - \frac{m}{L}}}}{1-\frac{m}{L}} \ge 1 +\sqrt{2\frac mL} + \frac mL \ge 1.\label{r}
    \end{equation}
\end{lemma}
The proof is essentially the same as that given in~\cite[Theorem~6]{U25}. For completeness, we give a proof in \cref{app:convergence}.

\begin{lemma}\label{lem:XX}
    Assume \cref{asm:C}, \cref{asm:F}, and \cref{asm:H} with $m \neq 0$. Then, for Wasserstein SR2FISTA~\eqref{SR2FISTA}, the following bound holds for any $k\ge1$:
    \[
        \normL{X_k - X^\star}{\nu}^2 \le \frac{4}{ma_k}W_2^2(\mu_0,\mu^\star).
    \]
\end{lemma}
The proof can be carried out by a straightforward extension of the Euclidean SR2FISTA case~\cite[Theorem~6]{U25} to the $L^2(\nu)$ case, using~\eqref{TWNAG1} with the Lyapunov function in~\cref{lem:lyapu2}.
For completeness, we give a proof in \cref{app:convergence}.

\begin{theorem}[Convergence rate of Wasserstein SR2FISTA]\label{thm:SR2FISTA}
    Assume \cref{asm:C}, \cref{asm:F}, and \cref{asm:H}. Then, for Wasserstein SR2FISTA~\eqref{SR2FISTA}, the following bound holds for any $k\ge1$:
    \[
        \GG(\mu_k) - \GG(\mu^\star) \le \begin{cases}
            \dfrac{2L}{k^2}W_2(\mu_0,\mu^\star)^2 & \text{if } m = 0,\\
            \min\setE*{\dfrac{6L}{k^2},\dfrac{3(L-m)}{2}r^{-(k-1)}} W_2(\mu_0,\mu^\star)^2 & \text{if } m > 0,
        \end{cases}
    \]
    where $r$ is defined by~\eqref{r}.
\end{theorem}
\begin{proof}
    From \cref{lem:lyapu2}, we have for any $k \ge 0$
    \[
        a_k\paren*{\GG(\mu_k) - \GG(\mu^\star) - \frac m2 \normL{X_k - X^\star}{\nu}^2} \le \mathcal E_k \le \mathcal E_0 = W_2^2(\mu_0,\mu^\star).
    \]
    When $m=0$, the claim follows directly from \cref{lem:SR2FISTAa}.
    When $m>0$, from~\cref{lem:XX}, we have
    \[
        \GG(\mu_k) - \GG(\mu^\star) \le \frac{1}{a_k}W_2^2(\mu_0,\mu^\star) + \frac m2 \normL{X_k - X^\star}{\nu}^2 \le \frac3{a_k} W_2^2(\mu_0,\mu^\star).
    \]
    Therefore, from \cref{lem:SR2FISTAa}, the claim follows.
\end{proof}

\section{Proof of \cref{prop:convex}}\label{sec:convex}
In this section, we first show that if $\FF$ is geodesically convex on $\CC$, $\tF$ is convex on a certain set of simple functions contained in $\CCnu$.
Then, we prove \cref{prop:convex} by passing to the limit.

In the following, a finite collection of Borel subsets of $\RR^d$, denoted by $\mathscr A=\setE{A_1,\ldots,A_N}$, is said to be a \textit{partition} if $\bigcup_{i=1}^N A_i=\RR^d$, $A_i\cap A_j=\varnothing$ for all $1\leq i<j\leq N$, and $\nu(A_i)>0$ for all $1\leq i\leq N$.

\begin{definition}
Let $\CC\subseteq\PP_2(\RR^d)$.
For a partition $\mathscr A=\setE{A_1,\ldots,A_N}$ and $X\in\CCnu$, define $X_{\mathscr A}\in L^2(\nu)$ by
\[
    X_{\mathscr A}\coloneqq\sum_{i=1}^N\paren*{\frac1{\nu(A_i)}\int_{A_i}X(x)\dd\nu(x)}\mathbf1_{A_i}.
\]
Let $\CCnu^{\mathrm d}\subseteq L^2(\nu)$ denote the set of all such $X_{\mathscr A}$ obtained from some partition $\mathscr A$ and some $X\in\CCnu$.
\end{definition}

\begin{lemma}\label{lem:CCd}
For $\CC \subseteq \PP_2(\RR^d)$, assume that~\cref{asm:C} holds. Then, we have $\CCnu^{\mathrm d}\subseteq\CCnu$. Moreover, $\CCnu^{\mathrm d}$ is dense in $\CCnu$ with respect to the $L^2(\nu)$ norm.
\end{lemma}

\begin{proof}
We show that $\CCnu^{\mathrm d}\subseteq\CCnu$.
Let $X_{\mathscr A}\in\CCnu^{\mathrm d}$, where $X\in\CCnu$ and $\mathscr A=\setE{A_1,\ldots,A_N}$.
For each $1\leq i\leq N$, let $\nu_i$ denote the normalized restriction of $\nu$ to $A_i$, defined by
\[
    \nu_i(B)\coloneqq\frac{\nu(B\cap A_i)}{\nu(A_i)}
\]
for any Borel set $B\subseteq\RR^d$.
Since the probability space $(A_i,\nu_i)$ is non-atomic, we can take independent $\RR^d$-valued random variables $(Z_i^{(k)})_{k\geq1}$ on this space such that $(Z_i^{(k)})_\#\nu_i=X_\#\nu_i$ for every $k\geq1$; see, e.g.,~\cite{K21b}.
For each $k\geq1$, define $Z^{(k)}$ on the whole reference space by $Z^{(k)}=Z_i^{(k)}$ on $A_i$.
Then, by $(Z_i^{(k)})_\#\nu_i=X_\#\nu_i$ and $X\in\CCnu$, we have
\[
    (Z^{(k)})_\#\nu=\sum_{i=1}^N\nu(A_i)((Z_i^{(k)})_\#\nu_i)=\sum_{i=1}^N\nu(A_i)(X_\#\nu_i)=X_\#\nu\in\CC.
\]
In particular, $Z^{(k)}\in L^2(\nu)$ because it has the same law as $X$.
Therefore, $Z^{(k)}\in\CCnu$.
Moreover, since $\CCnu$ is convex, for any integer $M\geq1$,
\[
    \frac1M\sum_{k=1}^M Z^{(k)}\in\CCnu.
\]

Since, for each $i$, the random variables $(Z_i^{(k)})_{k\geq1}$ are independent and identically distributed under $\nu_i$ with law $X_\#\nu_i$, we have
\begin{align}
    \MoveEqLeft
    \normL*{\frac1M\sum_{k=1}^M Z^{(k)}-X_{\mathscr A}}{\nu}^2\\
    &=\sum_{i=1}^N\nu(A_i)\int_{A_i}\norm*{\frac1M\sum_{k=1}^M Z_i^{(k)}(y)-\int_{A_i}X(x)\dd\nu_i(x)}^2\dd\nu_i(y)\\
    &=\sum_{i=1}^N\nu(A_i)\frac1M\int_{A_i}\norm*{X(y)-\int_{A_i}X(x)\dd\nu_i(x)}^2\dd\nu_i(y)\\
    &=\frac1M\normL{X-X_{\mathscr A}}{\nu}^2.
\end{align}
Thus, as $M\to\infty$, $\frac1M\sum_{k=1}^M Z^{(k)}\to X_{\mathscr A}$ in $L^2(\nu)$.
Since $\CCnu$ is closed, $X_{\mathscr A}\in\CCnu$.

We show that $\CCnu^{\mathrm d}$ is dense in $\CCnu$.
Let $X\in\CCnu$.
For each integer $N\geq1$, consider the $(2N^2)^d$ distinct cubes in $\RR^d$, denoted by $B_1,\ldots,B_{(2N^2)^d}$, of the form
\[
    \left[\frac{i_1}{N},\frac{i_1+1}{N}\right)\times\cdots\times\left[\frac{i_d}{N},\frac{i_d+1}{N}\right),\qquad (i_1,\ldots,i_d)\in\setE{-N^2,\ldots,N^2-1}^d.
\]
These cubes partition $[-N,N)^d$.
Let $B_0\coloneqq\RR^d\setminus[-N,N)^d$ and define
\[
    A_i\coloneqq\setI{x\in\RR^d}{X(x)\in B_i},\qquad i=0,\ldots,(2N^2)^d.
\]
If any of these sets have zero $\nu$-measure, merge all such sets into one of the positive-measure sets.
Let $\mathscr A_N$ denote the resulting partition.
% This modification affects the cells only by sets of $\nu$-measure zero.

For $\nu$-a.e. $y\in\RR^d\setminus A_0$, the vectors $X(y)$ and $X_{\mathscr A_N}(y)$ belong to the closure of the same cube.
Thus,
\[
    \norm{X(y)-X_{\mathscr A_N}(y)}\leq\frac{\sqrt d}{N}.
\]
Consequently,
\begin{align}
    \MoveEqLeft
    \normL{X-X_{\mathscr A_N}}{\nu}^2\\
    &=\int_{A_0}\norm{X-X_{\mathscr A_N}}^2\dd\nu+\int_{\RR^d\setminus A_0}\norm{X-X_{\mathscr A_N}}^2\dd\nu\\
    &\leq\int_{A_0}\norm{X-X_{\mathscr A_N}}^2\dd\nu+\frac d{N^2}\nu(\RR^d\setminus A_0).
\end{align}

If $\nu(A_0)>0$, set $m_0\coloneqq\nu(A_0)^{-1}\int_{A_0}X\dd\nu$.
Since $X_{\mathscr A_N}=m_0$ $\nu$-a.e. on $A_0$, we have
\[
    \int_{A_0}\norm{X-X_{\mathscr A_N}}^2\dd\nu=\int_{A_0}\norm{X}^2\dd\nu-\nu(A_0)\norm{m_0}^2\leq\int_{A_0}\norm{X}^2\dd\nu.
\]
The same inequality holds trivially when $\nu(A_0)=0$.
Since $A_0\subseteq\setI{x\in\RR^d}{\norm{X(x)}\geq N}$, it follows that
\[
    \normL{X-X_{\mathscr A_N}}{\nu}^2\leq\frac d{N^2}+\int_{\{\norm{X}\geq N\}}\norm{X}^2\dd\nu.
\]
The right-hand side converges to zero as $N\to\infty$ because $X\in L^2(\nu)$.
Therefore, $X_{\mathscr A_N}\to X$ in $L^2(\nu)$, proving that $\CCnu^{\mathrm d}$ is dense in $\CCnu$.
\end{proof}

\begin{lemma}\label{lem:geoconv0}
    For $\CC \subseteq \PP_2(\RR^d)$, assume \cref{asm:C}.
    For $\FF$ and $\tF$ defined by~\eqref{tF}, assume \cref{asm:F}.\ref{asm:diff}.
    Then, if $\FF$ is geodesically convex on $\CC$, $\tF$ is convex on $\CCnu^{\mathrm{d}}$.
\end{lemma}
The proof is similar to~\cite[Lemma~3.4]{P24}. The difference is that we use finite-valued maps in $\CCnu$, rather than in $L^2(\nu)$ in~\cite[Lemma~3.4]{P24}.
\begin{proof}
    Let $X_0,X_1\in\CCnu^\mathrm{d}$. By passing to a common refinement of their partition if necessary, we can write $X_0,X_1\in\CCnu^\mathrm{d}$ in the forms
\[
    X_0=\sum_{i=1}^N a_i\mathbf1_{A_i},\qquad X_1=\sum_{i=1}^N b_i\mathbf1_{A_i},\qquad p_i\coloneqq\nu(A_i)>0,
\]
for some $a_i ,b_i \in \RR^d$ ($i=1,\ldots,N$) and a common partition $\mathscr A = \setE{A_1,\ldots,A_N}$.
Set
\[
    x_i(t)\coloneqq (1-t)a_i+tb_i,\qquad X_t\coloneqq(1-t)X_0+tX_1,\qquad \rho_t\coloneqq(X_t)_\#\nu=\sum_{i=1}^Np_i\delta_{x_i(t)}.
\]
From \cref{lem:CCd}, we have $X_0,X_1 \in \CCnu$, and by convexity of $\CCnu$, we have $X_t\in\CCnu$. Hence,
\begin{equation}
    \rho_t\in\CC \quad  \text{for any }t\in[0,1]. \label{rhoinC}
\end{equation}
In the following, we show the convexity of $f(t) \coloneqq \FF(\rho_t)$ defined on $t \in [0,1]$, which yields the desired claim $\tF(X_t) \le (1-t)\tF(X_0) + t\tF(X_1)$ for any $t \in [0,1]$.

Fix $s\in[0,1]$ and let $\delta_s \coloneqq \min\setI{\norm{x_i(s) - x_j(s)}}{x_i(s) \neq x_j(s)}>0$. If $\rho_s$ has only one support point, set $\delta_s = +\infty$. 
Choose $h \in \RR$ so that $s+h \in [0,1]$ and 
\begin{equation}
    |h|\max_{1\leq i\leq N}\norm{b_i - a_i}<\frac{\delta_s}{2}. \label{hsmall}
\end{equation}
Then, the coupling
\[
    \gamma_{s,h}\coloneqq\sum_{i=1}^Np_i\delta_{(x_i(s),x_i(s+h))}
\]
is optimal between $\rho_s$ and $\rho_{s+h}$.
Indeed, if $\rho_s$ has only one support point, this is trivial.
Otherwise, each $x_i(s+h)$ is closer to $x_i(s)$ than to any other distinct support point $x_j(s)$ since 
\begin{align}
    \MoveEqLeft
    \norm{x_j(s) - x_i(s+h)} \ge \norm{x_j(s) - x_i(s)} - \norm{x_i(s) - x_i(s+h)}\\
    &\ge \delta_s - \norm{x_i(s) - x_i(s+h)}\\
    &> 2\abs{h}\norm{b_i - a_i} - \norm{x_i(s) - x_i(s+h)}\\
    &= \norm{x_i(s) - x_i(s+h)}.
\end{align}
Consequently, $\gamma_{s,h}$ minimizes the quadratic transport cost, and we have
\begin{equation}
    W_2^2(\rho_s,\rho_{s+h})=\sum_{i=1}^Np_i\norm{x_i(s+h) - x_i(s)}^2 = h^2\sum_{i=1}^Np_i\norm{a_i - b_i}^2.\label{localW}
\end{equation}
Moreover, when $h >0$, the geodesic from $\rho_s$ to $\rho_{s+h}$ is given by $\hat\rho_\tau = {T_\tau}_\# \gamma_{s,h}$ for $\tau \in [0,1]$, where $T_\tau(x,y) = (1-\tau)x + \tau y$.
By definition, we have $\hat\rho_\tau = \rho_{s + h\tau}$. Similarly, when $h < 0$, the geodesic from $\rho_{s+h}$ to $\rho_{s}$ coincides with $\rho_{s + h - h\tau}$.

Therefore, for every $s\in(0,1)$, there exists $\varepsilon>0$ such that the restrictions of $f$ to $[s-\varepsilon,s]$ and $[s,s+\varepsilon]$ are convex since $\FF$ is geodesically convex on $\CC$ and $\rho_t \in \CC$ for any $t \in [0,1]$ by~\eqref{rhoinC}.

For such $\varepsilon$, we next show that the restriction of $f$ to $[s - \varepsilon,s + \varepsilon]$ is convex.
To show this, it suffices to prove that $f$ is differentiable on $(0,1)$, since for differentiable $g \colon I \to \RR$ where $I \subseteq \RR$ is an open interval, $g$ is convex if and only if $g'$ is nondecreasing; see also~\cite[Lemma~3.8]{P24}.
For any $s \in (0,1)$, take $h \in \RR$ satisfying~\eqref{hsmall} again so that $\gamma_{s,h}$ is the optimal coupling.
% If $X_0=X_1$ in $L^2(\nu)$, then $\FF(\rho_t)$ is constant.
Wasserstein differentiability at $\rho_s$, applied to $\gamma_{s,h}$, gives
\begin{align}
    f(s+h) - f(s) &= \FF(\rho_{s+h})-\FF(\rho_s)\\
    &= \int_{\RR^d \times \RR^d} \inpr{\nablaw \FF(\rho_s)(x)}{y-x} \dd \gamma_{s,h}(x,y) + \mathrm{o}(W_2(\rho_{s+h},\rho_s))\\
    &=h\sum_{i=1}^Np_i\inpr{\nablaw\FF(\rho_s)(x_i(s))}{b_i - a_i} +o(|h|),
\end{align}
where~\eqref{localW} is used in the last equality.
Thus $f$ is differentiable on $(0,1)$.
Therefore, for any $s \in (0,1)$, there exists a neighborhood $U \subseteq [0,1]$ of $s$ such that the restriction of $f$ to $U$ is convex.
% \[
%     \frac{\dd}{\dd s}\FF(\rho_s) =\sum_{i=1}^Np_i\inpr{\nablaw\FF(\rho_s)(x_i(s))}{b_i - a_i}.
% \]

This, together with $f$ being continuous on $[0,1]$, implies that $f$ is convex. Therefore, the claim follows.
\end{proof}

\begin{lemma}\label{lem:geoconvc}
    For $\CC \subseteq \PP_2(\RR^d)$, assume \cref{asm:C}.
    For $\FF$ and $\tF$ defined by~\eqref{tF}, assume \cref{asm:F}.\ref{asm:diff}.
    Then, $\FF$ is geodesically convex on $\CC$ if and only if $\tF$ is convex on $\CCnu$.
\end{lemma}
\begin{proof}
Assume that $\FF$ is geodesically convex on $\CC$.
Let $X_0,X_1\in\CCnu$ be arbitrary.
From~\cref{lem:CCd}, for $ i = 0,1$, there exists a sequence $\setE{X_i^{(k)}}_{k\ge1} \subset \CCnu^{\mathrm{d}}$ such that 
\[
    X_i^{(k)}\to X_i\quad\text{in }L^2(\nu)\quad \text{as} \quad k \to \infty.
\]

For $t\in[0,1]$ and $k \ge 1$, define
\[
    X_t^{(k)}\coloneqq(1-t)X_0^{(k)}+tX_1^{(k)},\qquad X_t\coloneqq(1-t)X_0+tX_1.
\]
All these maps belong to $\CCnu$ by convexity, and
\[
    \normL{X_t^{(k)}-X_t}{\nu}\leq(1-t)\normL{X_0^{(k)}-X_0}{\nu}+t\normL{X_1^{(k)}-X_1}{\nu}\to 0.
\]
Therefore,
\[
    W_2((X_t^{(k)})_\#\nu,(X_t)_\#\nu)\leq\normL{X_t^{(k)}-X_t}{\nu}\to 0.
\]
Since $\FF$ is continuous at $(X_t)_\#\nu\in\CC\subseteq\CC'$,
\[
    \tF(X_t^{(k)})\to \tF(X_t).
\]
Applying~\cref{lem:geoconv0} to $X_0^{(k)},X_1^{(k)}$, we have
\[
    \tF((1-t)X_0^{(k)}+tX_1^{(k)})\leq(1-t)\tF(X_0^{(k)})+t\tF(X_1^{(k)}),
\]
and passing to the limit yields
\[
    \tF((1-t)X_0+tX_1)\leq(1-t)\tF(X_0)+t\tF(X_1).
\]
Thus $\tF$ is convex on $\CCnu$.

The converse is immediate.
\end{proof}

\begin{proof}[Proof of \cref{prop:convex}]
We assume $\FF$ is $m$-geodesically convex on $\CC$ and show that $\tF$ is $m$-convex on $\CCnu$.
For $\rho \in \PP_2(\RR^d)$, define
\[
    \mathcal M(\rho)\coloneqq\frac12\int\norm{x}^2\,d\rho(x),\qquad \RRc(\rho)\coloneqq\FF(\rho)-m\mathcal M(\rho).
\]

Let $\rho_0, \rho_1 \in \CC$ and $\gamma \in \Gamma_\mathrm{o}(\rho_0, \rho_1)$. For the Wasserstein geodesic $\setE{\rho_t}_{t\in[0,1]}$ from $\rho_0$ to $\rho_1$, we have
\begin{align}
    % \MoveEqLeft
    \mathcal M(\rho_t)&= \frac12 \int \norm{x}^2 \dd \rho_t(x)\\
    &= \frac12 \int \norm{(1-t)x + t y}^2 \dd \gamma(x,y)\\
    &= \frac12 \int \paren*{(1-t)\norm{x}^2 + t\norm{y}^2 - t(1-t)\norm{x - y}^2} \dd \gamma(x,y)\\
    &= (1-t)\mathcal M(\rho_0)+t\mathcal M(\rho_1)-\frac12t(1-t)W_2^2(\rho_0,\rho_1).
\end{align}
Since $\FF$ is $m$-geodesically convex on $\CC$,
$\RRc$ satisfies
\[
    \RRc(\rho_t) \le (1-t)\RRc(\rho_0) + t\RRc(\rho_1),
\]
implying $\RRc$ is geodesically convex on $\CC$.
Moreover, $\RRc$ is Wasserstein differentiable on $\CC'$.

Therefore, from \cref{lem:geoconvc}, the lift $\tR \colon L^2(\nu) \to \RR\cup\setE{+\infty}$ defined by
\[
    \tR(X) \coloneqq \RRc(X_\# \nu) = \tF(X)-\frac m2\normL{X}{\nu}^2
\]
is convex on $\CCnu$.

Thus, for $X_0,X_1 \in \CCnu$, we have
\[
    \tR((1-t)X_0+tX_1) \le (1-t)\tR(X_0) + t\tR(X_1).
\]
Using
\[
    (1-t)\normL{X_0}{\nu}^2+t\normL{X_1}{\nu}^2-\normL{(1-t)X_0+tX_1}{\nu}^2=t(1-t)\normL{X_0-X_1}{\nu}^2,
\]
we obtain
\[
    \tF((1-t)X_0+tX_1)\leq(1-t)\tF(X_0)+t\tF(X_1)-\frac m2t(1-t)\normL{X_0-X_1}{\nu}^2.
\]
Therefore, $\tF$ is $m$-convex on $\CCnu$.

Conversely, let $\mu_0,\mu_1\in\CC$ and $\gamma\in\Gamma_\mathrm{o}(\mu_0,\mu_1)$.
Since $\nu$ is non-atomic, there exist $X_0,X_1\in L^2(\nu)$ such that $(X_0,X_1)_\#\nu=\gamma$.
Then $X_0,X_1\in\CCnu$ and $\normL{X_0-X_1}{\nu}^2=W_2^2(\mu_0,\mu_1)$, so the desired inequality follows from the $m$-convexity of $\tF$.
\end{proof}

\section{Discussion and conclusion}

We developed Wasserstein FISTA and Wasserstein SR2FISTA for composite optimization over probability measures, establishing an $\Order{k^{-2}}$ objective-gap bound in the convex setting and accelerated linear convergence when the smooth component is geodesically strongly convex. These guarantees accommodate nonsmooth, extended-valued regularizers, including entropy. This extends discrete-time acceleration guarantees to problems not directly covered by existing analyses~\cite{T23,TLG26}.

A main limitation is that the proposed algorithms include a proximal map that generally has to be approximated numerically, adding computational cost. Our convergence guarantees assume that this proximal map is computed exactly, as in the analysis of Salim et al.~\cite{SKL20}. 
Understanding how approximation errors affect convergence remains an important problem. 
Another direction is to develop explicit methods that avoid such proximal computations while retaining discrete-time acceleration guarantees for such entropy-regularized problems. 
The insights gained in this work provide a starting point for developing and analyzing such methods.

\backmatter

\bmhead{Acknowledgements}
Kansei Ushiyama acknowledges OpenAI for providing complimentary access to ChatGPT through the ChatGPT for Academic Researchers program.

\section*{Declarations}

\bmhead{Competing interests}
The authors have no competing interests to declare that are relevant to the content of this article.

\bmhead{Use of LLMs}
The authors used a large language model (LLM) to assist with proving \cref{prop:convex} and~\eqref{meanentprox}, constructing examples in \cref{sec:ex}, and proofreading the English text of the manuscript.

% Some journals require declarations to be submitted in a standardised format. Please check the Instructions for Authors of the journal to which you are submitting to see if you need to complete this section. If yes, your manuscript must contain the following sections under the heading `Declarations':

% \begin{itemize}
% \item Conflict of interest/Competing interests (check journal-specific guidelines for which heading to use)
% \item Ethics approval and consent to participate
% \item Consent for publication
% \item Materials availability
% \end{itemize}

% \noindent
% If any of the sections are not relevant to your manuscript, please include the heading and write `Not applicable' for that section. 

\begin{appendices}

\section{Proof of~\eqref{meanentprox}}\label{app:meanconst}
Fix $Y \in L^2(\nu)$ and define
\[
    J_Y(X) \coloneqq \beta \ent(X_\# \nu)
    + \frac{1}{2\eta} \normL{X - Y}{\nu}^2,
    \qquad
    \bar X \coloneqq \int X \dd \nu,
    \qquad
    \bar Y \coloneqq \int Y \dd \nu.
\]
Then
\[
    \prox_{\eta\Phi}(Y)
    = \argmin_{X \in \CCnu \cap \KK_\nu} J_Y(X).
\]
To prove~\eqref{meanentprox}, it therefore suffices to show that
\[
    \argmin_{X \in \CCnu \cap \KK_\nu} J_Y(X)
    =
    \left(\argmin_{X \in \KK_\nu} J_Y(X)\right)
    + P_B(\bar Y) - \bar Y.
\]
We first establish two claims.

\begin{claim}
    For every $X \in L^2(\nu)$ and $a \in \RR^d$, we have
    \[
        \ent(X_\# \nu) = \ent((X + a)_\# \nu),
        \quad \text{and} \quad
        X \in \KK_\nu
        \Longleftrightarrow
        X + a \in \KK_\nu.
    \]
\end{claim}
\begin{proof}
    The entropy identity follows directly from the definition of $\ent$.
    To prove the second assertion, let $X \in \KK_\nu$ and define
    $\gamma_0 \coloneqq (\id,X)_\# \nu$,
        $\gamma_a \coloneqq (\id,X+a)_\# \nu.$
    For any $\gamma \in \Gamma(\nu,(X+a)_\# \nu)$, the shifted coupling
    $
        \widetilde{\gamma}
        \coloneqq (\id,\id-a)_\# \gamma
    $
    belongs to $\Gamma(\nu,X_\# \nu)$.
    Since $\gamma_0 \in \Gamma_{\mathrm{o}}(\nu,X_\# \nu)$, we obtain
    \begin{align}
        \MoveEqLeft
        \int \norm{x-y}^2 \dd \gamma_a(x,y)
        = \int \norm{x-(y+a)}^2 \dd \gamma_0(x,y) \\
        &=
        \int \norm{x-y}^2 \dd \gamma_0(x,y)
        - 2\inpr*{\int x \dd \nu(x) - \bar X}{a}
        + \norm{a}^2 \\
        &\le
        \int \norm{x-y}^2 \dd \widetilde{\gamma}(x,y)
        - 2\inpr*{\int x \dd \nu(x) - \bar X}{a}
        + \norm{a}^2 \\
        &=
        \int \norm{x-(y+a)}^2 \dd \widetilde{\gamma}(x,y) \\
        &=
        \int \norm{x-y}^2 \dd \gamma(x,y).
    \end{align}
    Since $\gamma$ was arbitrary, we conclude that
    $\gamma_a \in \Gamma_{\mathrm{o}}(\nu,(X+a)_\# \nu)$
    and hence $X+a \in \KK_\nu$.
    The converse follows by applying the same argument to $X+a$ and $-a$.
\end{proof}

\begin{claim}
    For $X \in \KK_\nu$, define
    \[
        U(X) \coloneqq X - \bar X + \bar Y.
    \]
    Then,
    \begin{equation}
        J_Y(X)
        = J_Y(U(X))
        + \frac{1}{2\eta}\norm{\bar X - \bar Y}^2.
        \label{JyX}
    \end{equation}
\end{claim}
\begin{proof}
    By the first claim,
    \[
        \ent(X_\# \nu) = \ent(U(X)_\# \nu).
    \]
    Moreover, the decomposition
    \begin{align}
        \normL{X-Y}{\nu}^2
        &=
        \normL{(X-\bar X)-(Y-\bar Y)}{\nu}^2
        + \norm{\bar X-\bar Y}^2 \\
        &=
        \normL{U(X)-Y}{\nu}^2
        + \norm{\bar X-\bar Y}^2
    \end{align}
    yields~\eqref{JyX}.
\end{proof}
\begin{proof}[Proof of~\eqref{meanentprox}]
    Let $X_0 \in \argmin_{X \in \KK_\nu} J_Y(X)$ and define
    \begin{equation}
        X_B \coloneqq X_0 + P_B(\bar Y) - \bar Y.\label{XB}
    \end{equation}
    First, we show that $X_B \in \argmin_{X \in \KK_\nu \cap \CCnu} J_Y(X)$.
    By the first claim, $X_B \in \KK_\nu$.
    By~\eqref{JyX} and the optimality of $X_0$, we have $\bar X_0 = \bar Y$, and hence
    \[
        \bar X_B
        = \bar X_0 + P_B(\bar Y) - \bar Y
        = P_B(\bar Y) \in B.
    \]
    Thus, $X_B \in \CCnu$, and by the definition of $X_B$,
    \[
        U(X_B)
        = X_B - \bar X_B + \bar Y
        = X_0.
    \]

    For any $X \in \KK_\nu \cap \CCnu$, we have
    $U(X) \in \KK_\nu$ by the first claim, and $\bar X \in B$.
    Therefore, by~\eqref{JyX}, the optimality of $X_0$,
    and the definition of $P_B$, we obtain
    \begin{align}
        J_Y(X)\label{XtoXB}
        &= J_Y(U(X))
        + \frac{1}{2\eta}\norm{\bar X - \bar Y}^2 \\
        &\ge J_Y(X_0)
        + \frac{1}{2\eta}\norm{\bar X - \bar Y}^2 \\
        &\ge J_Y(X_0)
        + \frac{1}{2\eta}\norm{P_B(\bar Y) - \bar Y}^2 \\
        &= J_Y(X_B).
    \end{align}
    This proves that
    $X_B \in \argmin_{X \in \KK_\nu \cap \CCnu} J_Y(X)$.

    Conversely, let
    $X_* \in \argmin_{X \in \KK_\nu \cap \CCnu} J_Y(X)$. We show that $X_* - P_B(\bar Y) + \bar Y \in \argmin_{X \in \KK_\nu} J_Y(X)$. 
    For any $X_0 \in \argmin_{X \in \KK_\nu} J_Y(X)$ and $X_B$ defined by~\eqref{XB}, $X_B \in \argmin_{X \in \KK_\nu \cap \CCnu} J_Y(X)$ and \eqref{XtoXB} holds for any $X \in \KK_\nu\cap\CCnu$.
    Thus, $J_Y(X_*) = J_Y(X_B)$, and substituting $X = X_*$ into~\eqref{XtoXB} shows that equality holds in both inequalities in~\eqref{XtoXB}.
    Therefore, we have $J_Y(U(X_*)) = J_Y(X_0)$ and $\bar X_* = P_B(\bar Y)$, where the latter follows from the uniqueness of the projection. These imply that
    \[
        U(X_*) \in \argmin_{X \in \KK_\nu} J_Y(X), \qquad U(X_*) = X_* - \bar X_* + \bar Y = X_* - P_B(\bar Y) + \bar Y.
    \]
    Combining these, we obtain $X_* - P_B(\bar Y) + \bar Y \in \argmin_{X \in \KK_\nu} J_Y(X)$.
\end{proof}

\section{Proofs of the results in \cref{sec:convergence}}\label{app:convergence}
\bmhead{Proof of \cref{lem:FISTAa}}
\begin{proof}
    We show $a_k \ge k^2/(4L)$ and $a_k \ge \paren*{1-\sqrt{\frac mL}}^{-k}/(2L)$.
    From~\eqref{FISTAa}, we have
    \[
        a_{k+1} = \frac{1+2La_k+\sqrt{1+4La_k(1+ma_k)}}{2(L-m)} \ge \frac{\frac12+2La_k+\sqrt{4La_k}}{2L} = \paren*{\frac{1}{2\sqrt L} + \sqrt{a_k}}^2,
    \]
    which, together with $a_0 = 0$, implies
    $
        \sqrt{a_k} \ge {k}/({2\sqrt{L}}),
    $
    i.e., $a_k \ge k^2/(4L)$.
    From~\eqref{FISTAa}, we have
    \[
        a_1 = \frac {1}{L-m} = \frac 1L \frac{1}{1-\frac mL} = \frac 1L \frac{1}{(1 + \sqrt{\frac mL})(1 - \sqrt{\frac mL})} \ge \frac {1}{2L(1 - \sqrt{\frac mL})},
    \]
    and
    \[
        a_{k+1} \ge \frac{2La_k+\sqrt{4Lma_k^2}}{2(L-m)} = \frac{1 + \sqrt{\frac mL}}{1-\frac mL}a_k = \frac{1}{1 - \sqrt{\frac mL}}a_k.
    \]
    These imply $a_k \ge \paren*{1-\sqrt{\frac mL}}^{-k}/(2L)$.
\end{proof}

\bmhead{Proof of \cref{lem:SR2FISTAa}}
\begin{proof}
    We show $a_k \ge (k/\sqrt{2L})^2$ and $a_k \ge \frac{2}{L - m}r^{k-1}$.
    First, by~\eqref{SR2FISTAa}, we have
    \[
        a_{k+1} =\frac{1+La_k+\sqrt{(1+ma_k)(1+(2L-m)a_k)}}{L-m} \ge \frac{\frac12+La_k+\sqrt{2La_k}}{L} =  \paren*{\sqrt{a_k} + \frac{1}{\sqrt{2L}}}^2.
    \]
    Since $a_0 = 0$, this implies the quadratic growth $a_k \ge (k/\sqrt{2L})^2$.
    Second, by~\eqref{SR2FISTAa}, we have $a_1 = \frac{2}{L - m}$ and
    \begin{align}
        &a_{k+1} \ge \frac{L + \sqrt{m(2L - m)}}{L - m}a_k = ra_k.
    \end{align}
    These imply $a_k \ge \frac{2}{L - m}r^{k-1}$.
\end{proof}

\bmhead{Proof of \cref{lem:XX}}
\begin{proof}    
    Let $U_k \coloneqq a_k(X_k - X^\star)$.
    Then, we have 
    \begin{align}
        U_{k+1}-U_k &= a_{k+1}(X_{k+1} - X^\star) - a_k(X_k-X^\star) \label{UU}\\
        &= (a_{k+1} - a_k)(X_{k+1}- X^\star) + a_k(X_{k+1} - X_k)\\
        &= (a_{k+1} - a_k)\paren*{X_{k+1}- X^\star + \frac{a_k}{a_{k+1} - a_k}(X_{k+1} - X_k)}\\
        & \overset{(*)}{=} (a_{k+1} - a_k)\paren*{X_{k+1}- X^\star + V_{k+1} - X_{k+1}}\\
        &= (a_{k+1} - a_k)(V_{k+1} - X^\star), 
    \end{align}
    where in the equality marked by $(*)$, we apply the update~\eqref{TWNAG1} to $X_{k+1} - X_k$.
    Since $\tF+\Phi$ is $m$-strongly convex and minimized at $X^\star$, we have $\GG(\mu_k)-\GG(\mu^\star)-\frac m2\normL{X_k-X^\star}{\nu}^2\ge0$. Thus, \cref{lem:lyapu2} implies that
    \[
        (1+ma_k)\normL{V_k - X^\star}{\nu}^2 \le \mathcal E_k \le \mathcal E_0 = \normL{\id - X^\star}{\nu}^2 = W_2^2(\mu_0,\mu^\star).
    \]
    Combining this with~\eqref{UU}, we have
    \begin{equation}
        \normL{U_{k+1} - U_k}{\nu} = (a_{k+1} - a_k)\normL{V_{k+1} -X^\star}{\nu} \le \frac{a_{k+1} - a_k}{\sqrt{1+ma_{k+1}}}W_2(\mu_0,\mu^\star).
    \end{equation}
    Here, the right-hand side can be bounded by using the inequality
    \[
        \frac{a_{k+1} - a_k}{\sqrt{1+ma_{k+1}}} \le \frac{2}{m}(\sqrt{1+ma_{k+1}} - \sqrt{1+ma_{k}}),
    \]
    which follows from
    \begin{align}
        &2(\sqrt{1+ma_{k+1}} - \sqrt{1+ma_{k}})\sqrt{1+ma_{k+1}} - m(a_{k+1} - a_k)\\
        ={}&2(1+ma_{k+1}) - 2\sqrt{1+ma_{k}}\sqrt{1+ma_{k+1}} - (1+ma_{k+1}) + (1+ma_{k})\\
        ={}&(\sqrt{1+ma_{k+1}} - \sqrt{1+ma_{k}})^2 \ge 0.
    \end{align}
    This allows for the telescoping sum
    \begin{align}
    \MoveEqLeft
    \sum_{s=0}^{k-1} \normL{U_{s+1} - U_s}{\nu} \le \sum_{s=0}^{k-1} \frac{2}{m}(\sqrt{1+ma_{s+1}} - \sqrt{1+ma_{s}})W_2(\mu_0,\mu^\star)\\
    &{}= \frac{2}{m}(\sqrt{1+ma_{k}} - \sqrt{1+ma_{0}})W_2(\mu_0,\mu^\star)\\
    &{}= \frac{2}{m}(\sqrt{1+ma_{k}} - 1)W_2(\mu_0,\mu^\star),
    \end{align}
    and, using $U_0 = 0$, we obtain
    \begin{align}
        \MoveEqLeft 
        a_k\normL{X_k - X^\star}{\nu} = \normL*{U_k}{\nu} = \normL*{\sum_{s=0}^{k-1} (U_{s+1} - U_s) + U_0}{\nu}\\
        &\le \sum_{s=0}^{k-1} \normL{U_{s+1} - U_s}{\nu}
        \le \frac{2}{m}(\sqrt{1+ma_{k}} - 1)W_2(\mu_0,\mu^\star).
    \end{align}
    Therefore, since 
    \[
        \frac{2}{ma_k}(\sqrt{1+ma_{k}} - 1) = \frac{2}{\sqrt{1+ma_{k}} + 1} \le \frac{2}{\sqrt{ma_{k}}},
    \]
    we have
    \[
        \normL{X_k - X^\star}{\nu}^2 \le \frac{4}{ma_k}W_2^2(\mu_0,\mu^\star).
    \]
\end{proof}




\end{appendices}


\bibliography{references_checked}% common bib file

@inproceedings{NEURIPS2025_347dae37,
 author = {Park, Jiyoung and Roy, Abhishek and Siegel, Jonathan W. and Bhattacharya, Anirban},
 booktitle = {Advances in Neural Information Processing Systems},
 doi = {10.52202/085713-1238},
 editor = {D. Belgrave and C. Zhang and H. Lin and R. Pascanu and P. Koniusz and M. Ghassemi and N. Chen},
 pages = {36828--36884},
 title = {Acceleration via silver step-size on Riemannian manifolds with applications to Wasserstein space},
 url = {https://proceedings.neurips.cc/paper\_files/paper/2025/file/347dae37f42ba41506cf7af753f98560-Paper-Conference.pdf},
 volume = {38, Main Conference},
 year = {2025}
}

@article {N13,
    AUTHOR = {Nesterov, Yu.},
     TITLE = {Gradient methods for minimizing composite functions},
   JOURNAL = {Math. Program.},
  FJOURNAL = {Mathematical Programming},
    VOLUME = {140},
      YEAR = {2013},
    NUMBER = {1, Ser. B},
     PAGES = {125--161},
      ISSN = {0025-5610},
   MRCLASS = {90C25 (68Q17 90C47)},
  MRNUMBER = {3071865},
MRREVIEWER = {Dana Simian},
       DOI = {10.1007/s10107-012-0629-5},
       URL = {https://doi.org/10.1007/s10107-012-0629-5},
}

@article {BT09,
    AUTHOR = {Beck, Amir and Teboulle, Marc},
     TITLE = {A fast iterative shrinkage-thresholding algorithm for linear
              inverse problems},
   JOURNAL = {SIAM J. Imaging Sci.},
  FJOURNAL = {SIAM Journal on Imaging Sciences},
    VOLUME = {2},
      YEAR = {2009},
    NUMBER = {1},
     PAGES = {183--202},
   MRCLASS = {35R30 (65F22 90C06 90C25)},
  MRNUMBER = {2486527},
       DOI = {10.1137/080716542},
       URL = {https://doi.org/10.1137/080716542},
}

@article{ADA21,
url = {http://dx.doi.org/10.1561/2400000036},
year = {2021},
volume = {5},
journal = {Foundations and Trends in Optimization},
title = {Acceleration Methods},
doi = {10.1561/2400000036},
issn = {2167-3888},
number = {1-2},
pages = {1--245},
author = {Alexandre d'Aspremont and Damien Scieur and Adrien Taylor}
}

@article {ADR22,
    AUTHOR = {Aujol, J.-F. and Dossal, Ch. and Rondepierre, A.},
     TITLE = {Convergence rates of the heavy ball method for quasi-strongly
              convex optimization},
   JOURNAL = {SIAM J. Optim.},
  FJOURNAL = {SIAM Journal on Optimization},
    VOLUME = {32},
      YEAR = {2022},
    NUMBER = {3},
     PAGES = {1817--1842},
      ISSN = {1052-6234,1095-7189},
   MRCLASS = {65K10 (90C25 90C30)},
  MRNUMBER = {4458909},
MRREVIEWER = {Driss\ Mentagui},
       DOI = {10.1137/21M1403990},
       URL = {https://doi.org/10.1137/21M1403990},
}

@ARTICLE{JGR25,
  title     = "Computer-assisted design of accelerated composite optimization
               methods: {OptISTA}",
  author    = {Jang, Uijeong and Das Gupta, Shuvomoy and Ryu, Ernest K.},
  journal   = {Math. Program.},
  publisher = "Springer Science and Business Media LLC",
  month     =  aug,
  year      =  2025,
  language  = "en",
  doi={10.1007/s10107-025-02258-5}
}

@article {CP16,
    AUTHOR = {Chambolle, Antonin and Pock, Thomas},
     TITLE = {An introduction to continuous optimization for imaging},
   JOURNAL = {Acta Numer.},
  FJOURNAL = {Acta Numerica},
    VOLUME = {25},
      YEAR = {2016},
     PAGES = {161--319},
      ISSN = {0962-4929,1474-0508},
   MRCLASS = {94A08 (49N90 65Kxx 68U05 68U10 90C25)},
  MRNUMBER = {3509209},
MRREVIEWER = {Jos\'e\ M.\ Guti\'errez},
       DOI = {10.1017/S096249291600009X},
       URL = {https://doi.org/10.1017/S096249291600009X},
}

@article {CC19,
    AUTHOR = {Calatroni, Luca and Chambolle, Antonin},
     TITLE = {Backtracking strategies for accelerated descent methods with
              smooth composite objectives},
   JOURNAL = {SIAM J. Optim.},
  FJOURNAL = {SIAM Journal on Optimization},
    VOLUME = {29},
      YEAR = {2019},
    NUMBER = {3},
     PAGES = {1772--1798},
      ISSN = {1052-6234,1095-7189},
   MRCLASS = {90C25 (65F22 94A12)},
  MRNUMBER = {3976651},
MRREVIEWER = {B\brevgrv ang\ C\^ong\ V\~u},
       DOI = {10.1137/17M1149390},
       URL = {https://doi.org/10.1137/17M1149390},
}

@ARTICLE{FV20,
  author={Florea, Mihai I. and Vorobyov, Sergiy A.},
  journal={IEEE Transactions on Signal Processing}, 
  title={A Generalized Accelerated Composite Gradient Method: Uniting {N}esterov's Fast Gradient Method and {FISTA}}, 
  year={2020},
  volume={68},
  pages={3033--3048},
  doi={10.1109/TSP.2020.2988614}}

@ARTICLE{FV19,
  author={Florea, Mihai I. and Vorobyov, Sergiy A.},
  journal={IEEE Transactions on Signal Processing}, 
  title={An Accelerated Composite Gradient Method for Large-Scale Composite Objective Problems}, 
  year={2019},
  volume={67},
  number={2},
  pages={444--459},
  doi={10.1109/TSP.2018.2866409}}

@article {P79,
    AUTHOR = {Passty, Gregory B.},
     TITLE = {Ergodic convergence to a zero of the sum of monotone operators
              in {H}ilbert space},
   JOURNAL = {J. Math. Anal. Appl.},
  FJOURNAL = {Journal of Mathematical Analysis and Applications},
    VOLUME = {72},
      YEAR = {1979},
    NUMBER = {2},
     PAGES = {383--390},
      ISSN = {0022-247X},
   MRCLASS = {47H05 (49A29)},
  MRNUMBER = {559375},
MRREVIEWER = {R.\ A.\ Plastock},
       DOI = {10.1016/0022-247X(79)90234-8},
       URL = {https://doi.org/10.1016/0022-247X(79)90234-8},
}

@article {RC22,
    AUTHOR = {Rebegoldi, Simone and Calatroni, Luca},
     TITLE = {Scaled, inexact, and adaptive generalized {FISTA} for strongly
              convex optimization},
   JOURNAL = {SIAM J. Optim.},
  FJOURNAL = {SIAM Journal on Optimization},
    VOLUME = {32},
      YEAR = {2022},
    NUMBER = {3},
     PAGES = {2428--2459},
      ISSN = {1052-6234,1095-7189},
   MRCLASS = {65K05 (90C25 90C30)},
  MRNUMBER = {4488588},
       DOI = {10.1137/21M1391699},
       URL = {https://doi.org/10.1137/21M1391699},
}

@article {DDD04,
    AUTHOR = {Daubechies, Ingrid and Defrise, Michel and De Mol, Christine},
     TITLE = {An iterative thresholding algorithm for linear inverse
              problems with a sparsity constraint},
   JOURNAL = {Comm. Pure Appl. Math.},
  FJOURNAL = {Communications on Pure and Applied Mathematics},
    VOLUME = {57},
      YEAR = {2004},
    NUMBER = {11},
     PAGES = {1413--1457},
      ISSN = {0010-3640,1097-0312},
   MRCLASS = {65J22 (35R30 47A52 65J20 65R30 90C30)},
  MRNUMBER = {2077704},
MRREVIEWER = {Hans-Peter\ Helfrich},
       DOI = {10.1002/cpa.20042},
       URL = {https://doi.org/10.1002/cpa.20042},
}

@article {RB16,
    AUTHOR = {Ryu, Ernest K. and Boyd, Stephen},
     TITLE = {A primer on monotone operator methods (survey)},
   JOURNAL = {Appl. Comput. Math.},
  FJOURNAL = {Applied and Computational Mathematics. An International
              Journal},
    VOLUME = {15},
      YEAR = {2016},
    NUMBER = {1},
     PAGES = {3--43},
      ISSN = {1683-3511,1683-6154},
   MRCLASS = {47H05 (47H09 47H10 47N10 65K05 90C25)},
  MRNUMBER = {3467167},
MRREVIEWER = {Radu\ Ioan\ Bo\c{t}},
}

@article{U25,
  title={A $\sqrt{2}$-accelerated {FISTA} for composite strongly convex problems},
  author={Ushiyama, Kansei},
  journal={arXiv preprint arXiv:2509.09295},
  year={2025}
}

@article {WL22,
    AUTHOR = {Wang, Yifei and Li, Wuchen},
     TITLE = {Accelerated information gradient flow},
   JOURNAL = {J. Sci. Comput.},
  FJOURNAL = {Journal of Scientific Computing},
    VOLUME = {90},
      YEAR = {2022},
    NUMBER = {1},
     PAGES = {11},
      ISSN = {0885-7474,1573-7691},
   MRCLASS = {65C05 (62F15 65M32)},
  MRNUMBER = {4342674},
MRREVIEWER = {Pavel\ T.\ Stoynov},
       DOI = {10.1007/s10915-021-01709-3},
       URL = {https://doi.org/10.1007/s10915-021-01709-3},
}

@article{CLTW25,
  author  = {Shi Chen and Qin Li and Oliver Tse and Stephen J. Wright},
  title   = {Accelerating optimization over the space of probability measures},
  journal = {Journal of Machine Learning Research},
  year    = {2025},
  volume  = {26},
  number  = {31},
  pages   = {1--40},
  url     = {http://jmlr.org/papers/v26/23-1288.html}
}

@article{TLG26,
  title={{Nesterov} acceleration in optimizing over probability measures},
  author={Tang, Jiaqi and Li, Qin and Gangbo, Wilfrid},
  journal={arXiv preprint arXiv:2607.23008},
  year={2026}
}

@inproceedings{MKLGSB21,
 author = {Mokrov, Petr and Korotin, Alexander and Li, Lingxiao and Genevay, Aude and Solomon, Justin and Burnaev, Evgeny},
 booktitle = {Advances in Neural Information Processing Systems},
 pages = {15243--15256},
 publisher = {Curran Associates, Inc.},
 title = {Large-Scale {Wasserstein} Gradient Flows},
 volume = {34},
 year = {2021}
}

@article {P24,
    AUTHOR = {Parker, Guy},
     TITLE = {Some convexity criteria for differentiable functions on the
              2-{W}asserstein space},
   JOURNAL = {Bull. Lond. Math. Soc.},
  FJOURNAL = {Bulletin of the London Mathematical Society},
    VOLUME = {56},
      YEAR = {2024},
    NUMBER = {5},
     PAGES = {1839--1858},
      ISSN = {0024-6093,1469-2120},
   MRCLASS = {49Q22 (26B25 39B62)},
  MRNUMBER = {4755333},
       DOI = {10.1112/blms.13030},
       URL = {https://doi.org/10.1112/blms.13030},
}

@article{GT19,
  author  = {Gangbo, Wilfrid and Tudorascu, Adrian},
  title   = {On differentiability in the {Wasserstein} space and well-posedness for {Hamilton--Jacobi} equations},
  journal = {Journal de Math{\'e}matiques Pures et Appliqu{\'e}es},
  volume  = {125},
  pages   = {119--174},
  year    = {2019},
  doi     = {10.1016/j.matpur.2018.09.003}
}

@inproceedings{SKL20,
  author = {Salim, Adil and Korba, Anna and Luise, Giulia},
  title = {The {Wasserstein} Proximal Gradient Algorithm},
  booktitle = {Advances in Neural Information Processing Systems},
  volume = {33},
  year = {2020},
}

@book {K21b,
    AUTHOR = {Kallenberg, Olav},
     TITLE = {Foundations of modern probability},
    SERIES = {Probability Theory and Stochastic Modelling},
    VOLUME = {99},
   EDITION = {Third},
 PUBLISHER = {Springer, Cham},
      YEAR = {2021},
      ISBN = {978-3-030-61871-1; 978-3-030-61870-4},
   MRCLASS = {60-01 (60A10 60G05)},
  MRNUMBER = {4226142},
MRREVIEWER = {Myron\ Hlynka},
       DOI = {10.1007/978-3-030-61871-1},
       URL = {https://doi.org/10.1007/978-3-030-61871-1},
}

@article{JKO98,
  author  = {Jordan, Richard and Kinderlehrer, David and Otto, Felix},
  title   = {The Variational Formulation of the {Fokker--Planck} Equation},
  journal = {SIAM Journal on Mathematical Analysis},
  volume  = {29},
  number  = {1},
  pages   = {1--17},
  year    = {1998},
  doi     = {10.1137/S0036141096303359},
  url     = {https://epubs.siam.org/doi/10.1137/S0036141096303359}
}

@inproceedings{CB18,
  author    = {Chizat, L{\'e}na{\"i}c and Bach, Francis},
  title     = {On the Global Convergence of Gradient Descent for Over-parameterized Models using Optimal Transport},
  booktitle = {Advances in Neural Information Processing Systems},
  volume    = {31},
  year      = {2018},
}

@book{AGS08,
  author    = {Ambrosio, Luigi and Gigli, Nicola and Savar{\'e}, Giuseppe},
  title     = {Gradient Flows: In Metric Spaces and in the Space of Probability Measures},
  edition   = {2},
  series    = {Lectures in Mathematics ETH Z{\"u}rich},
  publisher = {Birkh{\"a}user},
  year      = {2008},
  doi       = {10.1007/978-3-7643-8722-8},
  url       = {https://link.springer.com/book/10.1007/978-3-7643-8722-8}
}

@article{YCY24,
  author  = {Yao, Rentian and Chen, Xiaohui and Yang, Yun},
  title   = {{Wasserstein} Proximal Coordinate Gradient Algorithms},
  journal = {Journal of Machine Learning Research},
  volume  = {25},
  number  = {269},
  pages   = {1--66},
  year    = {2024},
  url     = {https://jmlr.org/papers/v25/23-0889.html}
}

@inproceedings{DBCS23,
  author    = {Diao, Michael Ziyang and Balasubramanian, Krishna and Chewi, Sinho and Salim, Adil},
  title     = {Forward-Backward {Gaussian} Variational Inference via {JKO} in the {Bures--Wasserstein} Space},
  booktitle = {Proceedings of the 40th International Conference on Machine Learning},
  series    = {Proceedings of Machine Learning Research},
  volume    = {202},
  pages     = {7960--7991},
  publisher = {PMLR},
  year      = {2023},
  url       = {https://proceedings.mlr.press/v202/diao23a.html}
}

@inproceedings{BUDAK24,
  author    = {Bonet, Cl{\'e}ment and Uscidda, Th{\'e}o and David, Adam and Aubin-Frankowski, Pierre-Cyril and Korba, Anna},
  title     = {Mirror and Preconditioned Gradient Descent in {Wasserstein} Space},
  booktitle = {Advances in Neural Information Processing Systems},
  volume    = {37},
  pages     = {25311--25374},
  year      = {2024},
  doi       = {10.52202/079017-0798},
}

@misc{T23,
  author        = {Tanaka, Ken'ichiro},
  title         = {Accelerated Gradient Descent Method for Functionals of Probability Measures by New Convexity and Smoothness Based on Transport Maps},
  year          = {2023},
  eprint        = {2305.05127},
  archivePrefix = {arXiv},
  primaryClass  = {math.OC},
  url           = {https://arxiv.org/abs/2305.05127}
}

@article{UTDT26,
  title={An optimal first-order method for smooth and strongly convex composite optimization and its stationary limit},
  author={Upadhyaya, Manu and Thomsen, Daniel Berg and Dieuleveut, Aymeric and Taylor, Adrien B},
  journal={arXiv preprint arXiv:2605.22929},
  year={2026}
}

@book{CD18,
  author    = {Carmona, Ren{\'e} and Delarue, Fran{\c{c}}ois},
  title     = {Probabilistic Theory of Mean Field Games with Applications {I}:
               Mean Field {FBSDEs}, Control, and Games},
  series    = {Probability Theory and Stochastic Modelling},
  volume    = {83},
  publisher = {Springer},
  address   = {Cham},
  year      = {2018},
  doi       = {10.1007/978-3-319-58920-6}
}

@misc{C13,
  author = {Cardaliaguet, Pierre},
  title  = {Notes on Mean Field Games},
  note   = {From {P.-L. Lions}' lectures at {Coll\`ege de France}},
  year   = {2013},
    url = {https://www.ceremade.dauphine.fr/%7Ecardaliaguet/MFG20130420.pdf}
}

@article{LMSS26,
  author  = {Lascu, Razvan-Andrei and Majka, Mateusz B. and
             {\v{S}}i{\v{s}}ka, David and Szpruch, {\L}ukasz},
  title   = {Linear Convergence of Proximal Descent Schemes
             on the Wasserstein Space},
  journal = {Annals of Applied Probability},
  year    = {2026}
}

\end{document}